%% file: agt-4-45.tex
\documentclass{gtart_h}
\input agtout

\lognumber{45}
\volumenumber{4}
\volumeyear{2004}
\papernumber{45}
\pagenumbers{1045}{1081}
\received{26 May 2003} 
\revised{7 October 2004}
\accepted{13 October 2004}
\published{6 November 2004}

\usepackage{amsmath,amssymb,graphicx}

\newtheorem{prop}{Proposition}
\newtheorem{theorem}{Theorem}

\newtheorem{corollary}{Corollary}
\newcommand{\oplusop}[1]{{\mathop{\oplus}\limits_{#1}}}

\def\sh#1{{\bf #1}\qua}

\def\C{\mathbb C}
\def\R{\mathbb R}
\def\Z{\mathbb Z}

\def\lra{\longrightarrow}
\def\Hom{\mathrm{Hom}}
\def\mc{\mathcal} 
\def\mf{\mathfrak}
\def\foams{\mathrm{Foams}}  
\def\bracket#1{\langle #1 \rangle}        
\newcommand{\bracketG}{\bracket{\Gamma}}   
\newcommand{\define}{\stackrel{\mbox{\scriptsize{def}}}{=}}
\def\drawing#1{\begin{center} \includegraphics{#1} \end{center}}
\def\drawings#1#2{\begin{center} \includegraphics[width=#2]{#1} \end{center}}
\def\cH{\mc{H}}     
\def\cA{\mc{A}}     
\def\cF{\mc{F}}     

\begin{document}
\title[sl(3) link homology]{sl(3) link homology} 
\author{Mikhail Khovanov} 
\address{Department of Mathematics, University of California, Davis CA 95616,
USA}
\email{mikhail@math.ucdavis.edu} 
 \begin{abstract} 
  We define a bigraded homology theory whose Euler 
  characteristic is the quantum sl(3) link invariant. 
 \end{abstract} 


\maketitle


\section{Introduction}

To each simple Lie algebra $\mf{g}$ there is associated a polynomial 
invariant $\bracket{L}$ of oriented links in $\R^3$ with components 
labelled by irreducible finite-dimensional 
representation of $\mf{g},$ see \cite{RT2}, \cite{Le} and references 
 therein. Suitably normalized, the invariant 
takes values in $\Z[q,q^{-1}].$  

We conjecture that, for each simply-laced $\mf{g},$ the 
invariant $\bracket{L}$ admits a categorification: there exists 
a bigraded homology theory of labelled links 
\begin{equation*}
\mc{H}_{\mf{g}}(L)= \oplus_{i,j\in \Z} \mc{H}_{\mf{g}}^{i,j}(L)
\end{equation*} 
with $\bracket{L}$ as the Euler characteristic, 
\begin{equation*} 
 \bracket{L} = \sum_{i,j}(-1)^i q^j \mbox{rk}
   (\mc{H}_{\mf{g}}^{i,j}(L)). 
\end{equation*}
In the first nontrivial case---when $\mf{g}$ is $\mf{sl}(2)$ and 
all components of $L$ are labelled by its defining representation---this 
theory was constructed in \cite{me:jones}. An oriented 
surface cobordism $S\subset \R^4$ between links $L_1$ and $L_2$ induces 
a homomorphism $\mc{H}(S): \mc{H}(L_1)\to \mc{H}(L_2),$ well-defined 
up to overall minus sign, see \cite{MagnusJ}, \cite{me:cobordisms}. 
These homomorphisms fit together 
into a functor from the category of link cobordisms to the category 
of bigraded abelian groups and group homomorphisms (with identifications 
 $f=-f$ for each homomorphism). Likewise, the theory 
$\mc{H}_{\mf{g}}$ should extend to a functor from the category of 
oriented framed cobordisms between framed links to the category of 
bigraded abelian groups and group homomorphisms, with $\pm$ identification. 
For $\mf{sl}(2)$ and its defining representation the framing can be 
hidden, but not, apparently, in the general case. Components of links 
and link cobordisms should be labelled by irreducible 
representations of $\mf{g}.$

The homology of the unknot labelled by an irreducible 
representation $V$ 
should be a commutative Frobenius algebra of dimension $\mbox{dim}(V),$
since the category of link cobordisms contains a subcategory equivalent 
to the category of 2-cobordisms between 1-manifolds 
(not embedded anywhere).  This algebra should be, in 
addition, graded, and its graded dimension should equal the quantum 
dimension of $V.$ When $V$ is the $k$-th exterior power of 
the defining representation of $\mf{sl}(n),$ this algebra is going to 
be the 
cohomology of the Grassmannian of $k$-dimensional subspaces of $\C^n.$ 
 
When $\mf{g}=\mf{sl}(2)$ and components of $L$ are labelled by arbitrary 
$\mf{sl}(2)$ representations, the invariant $\bracket{L}$ is the 
colored Jones polynomial. Its categorification is sketched 
in \cite{me:color}.

\begin{figure}[ht!]  \drawings{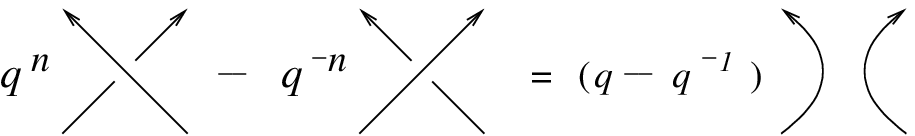}{8cm}\caption{Quantum 
$\mf{sl}(n)$ skein formula} \label{homfly} 
 \end{figure}

When $\mf{g}=\mf{sl}(n)$ and each components of $L$ is labelled either by 
the defining representation $V$ or its dual, the invariant is determined 
by the skein relation in Figure~\ref{homfly}. 
If we introduce a second variable $p=q^n$, the skein relation gives 
rise to the HOMFLY polynomial, a 2-variable polynomial invariant of 
oriented links \cite{Homfly}. 
We do not believe in a triply-graded homology 
theory categorifying the HOMFLY polynomial. Instead, for each $n\ge 0$ 
there should exist a bigraded theory categorifying the $(q,q^n)$ 
specialization of HOMFLY. For $n=0$ such theory was constructed 
by Peter Ozsv\'ath and Zolt\'an Szab\'o \cite{OZ1} and, independently, 
by Jacob Rasmussen \cite{R}. For $n=1$ the polynomial invariant as well as 
the homology theory is trivial. For $n=2$ the invariant is the Jones 
polynomial and was categorified in \cite{me:jones}. In this paper we 
construct link homology in the $n=3$ case. In the sequel we will 
extend this invariant to tangles and tangle cobordisms.

We use Greg Kuperberg's approach to  the quantum $\mf{sl}(3)$ link 
invariant \cite{Kspiders}. He described 
the representation theory of  $U_q(\mf{g}),$ for each rank two  
simple Lie algebra $\mf{g},$ via a calculus of 
planar trivalent graphs, called \emph{webs.} Each edge of a web is 
labelled by one of the two 
fundamental representation of $U_q(\mf{g}),$ 
and a trivalent vertex denotes a vector in the invariant space of a 
triple tensor product. 
A web with boundary gives a vector in the invariant 
space of a tensor product of fundamental representations. 
A closed web evaluates to a polynomial in $q^{\pm \frac{1}{2}}$ 
with integer coefficients if $\mf{g}=\mf{sl}(3);$ 
to a polynomial in  $q^{\pm 1}$ with integer 
coefficients, if $\mf{g}=\mf{so}(5);$ and to a 
rational function of $q$ if $\mf{g}=G_2$ (it is not immediately 
clear from \cite{Kspiders} that $G_2$ webs always evaluate to 
polynomials).  

Among $\mf{sl}(3), \mf{so}(5),$ and $G_2,$ 
the three simple Lie algebras or rank $2,$ only $\mf{sl}(3)$ is simply-laced. 
This property seems to relate to positivity of the 
Kuperberg calculus for $\mf{sl}(3).$ 
Namely, in the $\mf{sl}(3)$ case, closed webs
evaluate to Laurent polynomials in $q$ with \emph{non-negative} coefficients, 
after changing Kuperberg's $q^{\frac{1}{2}}$ to $-q.$ No such 
positivity is at hand for the other two calculi.

To compute the $\mf{sl}(3)$ invariant of a link $L$ choose its  
 plane projection $D$ and flatten each crossing in two possible ways, 
 as in Figure~\ref{flats}. This action produces $2^n$ closed webs, 
 where $n$ is the number of crossings of $D.$ 
 The invariant of $L$ is the sum of  these webs'  
 evaluations weighted by powers of $q$ and $-1,$ 
 as explained in Figure~\ref{cross}. Webs are evaluated recursively, 
 through Figure~\ref{rules} reductions (where 
 $[3]=q^2 + 1 + q^{-2}$ and 
 $[2]= q+ q^{-1}$). Each transformation either reduces the number of 
vertices or removes a loop. Since the formulas have no minus signs in 
them, the evaluation $\bracket{\Gamma}$ of any closed web $\Gamma$ has 
non-negative integer coefficients and lies in $\Z_+[q,q^{-1}].$ 

Positivity of $\bracketG$ is a hint that its categorification is just a graded 
abelian group, rather than a complex of graded abelian 
groups. Since $\bracketG= \sum a_j q^j$ for some non-negative integers 
$a_j,$ the categorification (or homology) $\cF(\Gamma)$ of 
$\bracketG$ should be a graded abelian group $\oplusop{j\in \Z} F_j$ 
 with $\mathrm{rank}(F_j) = a_j.$

\begin{figure}[ht!]  \drawings{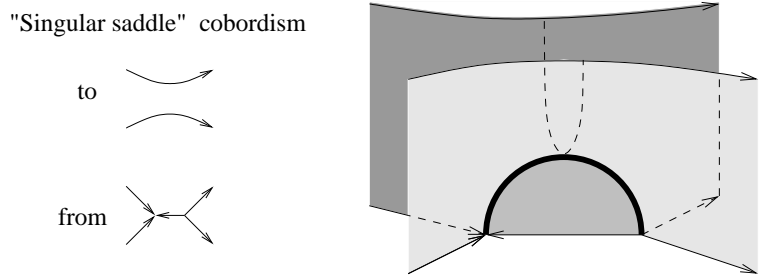}{10cm}\caption{Only the nontrivial  
 portion of the cobordism is indicated} 
\label{basic} \end{figure}

Homology $\cF(\Gamma)$ for different $\Gamma$ must be related. 
Namely, we need to have maps between homologies of any two 
webs that differ as in the R.H.S. of Figure~\ref{cross}. 
With these maps at hand we could set up a commutative cube 
of web homologies and maps, form the total complex of this cube 
and take its homology, following the approach of \cite{me:jones}. 
How to define these maps? In the $\mf{sl}(2)$ homology theory 
of links maps come from cobordisms. We want to have maps between 
homology of webs, and expect those to come from web cobordisms, 
which we will call foams. Locally, a foam should be homeomorphic to a 
neighbourhood of a point on the web times an interval. Since 
webs have singular points (vertices) with neighbourhoods
homeomorphic to letter $Y,$ foams will have singular arcs near 
which, locally, the foam is $Y\times [0,1].$ An example of a 
foam is depicted in Figure~\ref{basic}. Another example is the 
identity cobordism $\Gamma\times [0,1].$ 

To each foam $U$ with boundaries $\Gamma_1$ and $\Gamma_2$ we want 
to associate a homomorphism $\cF(U):\cF(\Gamma_1) \lra \cF(\Gamma_2).$

For the moment let us put foams aside and return to webs. 
The unknot evaluates to $[3]=q^2 + 1 + q^{-2}.$ 
Therefore, to the unknot we should associate a graded Frobenius algebra
 $\cA$ which has $\Z$ in degrees $-2, 0, 2.$ 
Our only choice is to take $\cA$ to be the ring $\Z[X]/X^3,$ the 
cohomology ring of $\mathbb{P}^2.$ We shift the degrees of 
this cohomology ring down by $2$ to balance them around degree $0.$

\begin{figure}[ht!]  \drawing{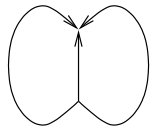}\caption{Theta web $\Theta$} 
\label{theta1} \end{figure}

\begin{figure}[ht!]  \drawing{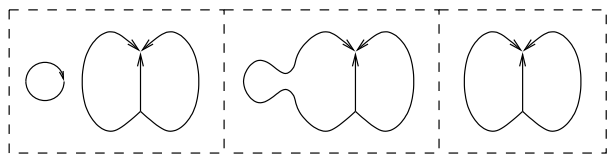}\caption{Web cobordism 
 cross-sections: a circle merges into an edge of $\Theta$} 
\label{theta2} \end{figure}
The theta web in Figure~\ref{theta1} is the simplest web with 
trivalent vertices. It evaluates to $[3][2],$ and we want to 
assign to $\Theta$ a graded abelian group $\cF(\Theta)$ 
of graded rank $[3][2].$ 
For each edge of $\Theta$ there is a web cobordism of merging 
a circle into the edge, depicted in Figure~\ref{theta2}. This cobordism 
should induce a map of abelian groups 
 \begin{equation*} 
   \cA \otimes \cF(\Theta) \lra \cF(\Theta) 
 \end{equation*} 
 and make $\cF(\Theta)$ into an $\cA$-module. The three 
 $\cA$-module structures on $\cF(\Theta)$ should commute and 
 make $\cF(\Theta)$ into an $\cA^{\otimes 3}$-module. 

\begin{figure}[ht!] \drawing{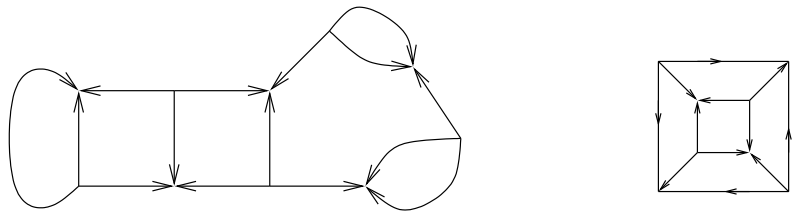}\caption{A digon web and 
 the cube web} 
\label{theta3} \end{figure}

Consider the cohomology ring of the full flag variety of $\C^3,$ 
 $$Fl_3\define
 \{ (V_1, V_2)| V_1\subset V_2 \subset \C^3, \mbox{dim}V_i=i \}.$$
Choose a hermitian inner product on $\C^3.$ The flag variety is 
homeomorphic to the space of triples of pairwise orthogonal complex 
lines in $\C^3$: 
$$Fl_3 \cong \{ (W_1, W_2, W_3) |W_i\subset \C^3, 
   \mbox{dim}W_i=1, \hspace{0.1in} W_i \bot W_j, i\not= j \}. $$ 
This homeomorphism defines an embedding of $Fl_3$ into 
$\mathbb{P}^2 \times \mathbb{P}^2 \times \mathbb{P}^2.$ The induced 
map on cohomology makes $\mbox{H}^{\ast}(Fl_3,\Z)$ into an 
$\cA^{\otimes 3}$-module. Note that $q^3 [3][2]$ is the graded dimension of 
 $\mbox{H}^{\ast}(Fl_3,\Z).$ We define $\cF(\Theta)$ as 
 the integral cohomology ring of $Fl_3$ with the grading shifted down 
by $3$: 
\begin{equation*} 
 \cF(\Theta) = \mathrm{H}^{\ast}(Fl_3,\Z)\{ -3\}.
\end{equation*} 
A naive generalization of this approach to $\cF(\Gamma)$ 
works for \emph{digon} webs, that is, webs which can be obtained from 
 a circle by adding digon faces, see Section~\ref{flagvar}.
 A digon web is shown in Figure~\ref{theta3} left. 
 To define $\cF(\Gamma)$ already for the "cube" web on 
 Figure~\ref{theta3} requires a different method. 

 In Section~\ref{sfoams} we define homology $\cF(\Gamma)$ for 
 any web $\Gamma.$ A foam $U$ should induce a map between 
 homology groups of its boundaries. The homology of the empty 
 web should be $\Z,$ and if $U$ has no boundaries (we call $U$ 
 a closed foam), it should induce an endomorphism of $\Z,$ that is, 
 an integer. These integers (evaluations of closed foams) must 
 be compatible with our assignments of $\cA$ to the circle web 
 and cohomology of the flag variety  to the theta web. Compatibility 
 requirements force us into a unique way to evaluate closed foams, 
 explained in Section~\ref{pre-foams}.  

 A foam which is a cobordism from the empty web to $\Gamma$ 
 should define a homomorphism from $\Z$ to $\cF(\Gamma),$ 
 or, equivalently, an element of $\cF(\Gamma).$ 
 We define homology $\cF(\Gamma)$ implicitly, as abelian 
 group generated by foams from the empty web to $\Gamma$ and  
 relations coming from required compatibility with 
 evaluations of closed foams, see Section~\ref{foams}.    
 The price of our implicit definition is in the work needed to 
 establish skein relations for web homology. This is done in 
 Section~\ref{webskein}.

Given a diagram $D$ of a link $L,$ we arrange the graded abelian groups 
$\cF(D_f),$ for all flattenings $f$ of $D,$ into a complex $\cF(D)$ 
of graded abelian groups, see Section~\ref{chains}. In 
Section~\ref{invariance} we prove 

\begin{theorem} If diagrams $D_1$ and $D_2$ are related by 
a Reidemeister move, the complexes $\cF(D_1)$ and $\cF(D_2)$ are 
chain homotopic. 
\end{theorem} 
Let $\cH(D) = \oplusop{i,j\in \Z}\cH^{i,j}(D)$ be the cohomology groups 
 of $\cF(D).$ 

\begin{corollary} Isomorphism classes of groups $\cH^{i,j}(D)$ are 
invariants of $L$.
\end{corollary} 

 \begin{prop} The Euler characteristic of $\cH(L)$ is the quantum 
 $\mf{sl}_3$ invariant: 
 \begin{equation*} 
   \bracket{L} = \sum_{i,j} (-1)^i q^j \mathrm{rk}(\cH^{i,j}(L)). 
 \end{equation*} 
 \end{prop} 

 Let us briefly discuss whether the invariants defined in 
 \cite{me:jones} and the present paper should be called \emph{homology} 
 or \emph{cohomology}. 
 We are building functors from the category of link cobordisms to the 
 category of abelian groups. Homology is a covariant functor, 
 cohomology a contravariant one. A link cobordism between two 
 links can be viewed as a morphism in either direction. 
 In other words, the category of link cobordisms has a 
 contravariant involution, which restricts to the identity on 
 objects. Composing with this involution turns any covariant 
 functor into a contravariant one. This allows 
 us to turn any homology theory of cobordisms into a cohomology 
 theory and vice versa. Previously we used the term "cohomology," 
 and in this paper switch to "homology," for balance. 
 Differentials in our complexes increase the index by $1,$ but 
 this convention can be easily reversed.

 A calculus of planar trivalent graphs that 
 gives rise to quantum $\mf{sl}(n)$ link invariants for any $n$ 
 was developed by Murakami, Ohtsuki and Yamada \cite{MOY}, 
 and further studied by Dongseok Kim \cite{Dongseok}. This 
 calculus has the same positivity flavor as Kuperberg's. 
 We expect that to each such planar graph there is naturally associated 
 a graded abelian group whose graded rank equals the polynomial 
 invariant of the graph, and these groups can be assembled 
 into complexes to produce quantum $\mf{sl}(n)$ link 
 homology. Cohomology rings of Grassmannians and partial flag 
 varieties will play an essential role in this theory.

\medskip 

{\bf Acknowledgments}\qua I am indebted to Greg Kuperberg 
 for many illuminating discussions and valuable input. 
 Section~\ref{flagvar} was written jointly with Greg Kuperberg.   
 Dror Bar-Natan's suggestion to try making this paper readable 
 led to a rather long introduction. While writing this paper, 
 I was partially supported by the NSF grant DMS-0104139 
 and daily intakes of grande white chocolate mocha at Starbucks.  
  

\section{Planar trivalent graphs and the quantum sl(3) invariant} 

\begin{figure} [ht!] \drawing{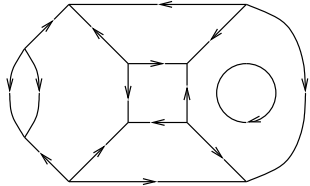}\caption{A closed web} 
 \label{web} 
 \end{figure}

A closed web in the Kuperberg $\mf{sl}(3)$ spider is a trivalent 
oriented graph in $\R^2,$ possibly 
with verticeless loops, such that at each vertex the edges are either all 
oriented in or oriented out, see Figure~\ref{web}. Any such graph  
is bipartite. 

\begin{figure} [ht!] \drawings{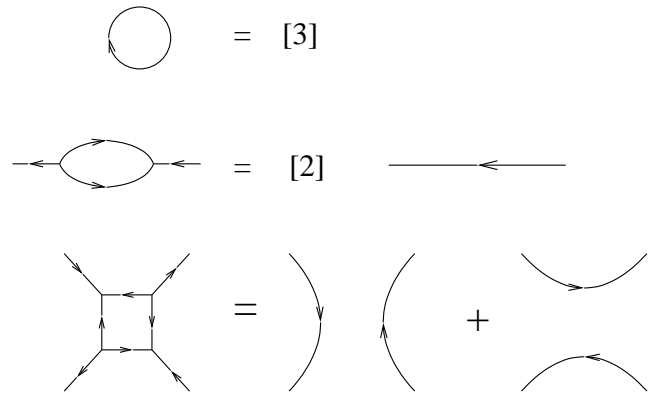}{8.5cm}\caption{Web skein relations} 
  \label{rules} 
 \end{figure}

To a web $\Gamma$ associate a Laurent polynomial 
$\bracket{\Gamma}\in \Z[q,q^{-1}],$ the \emph{Kuperberg bracket} of $\Gamma,$ 
via the rules depicted in Figure~\ref{rules}, where $[2]= q + q^{-1}$ and 
$ [3]= q^2+ 1 + q^{-2}.$ 

These rules say that 

(i)\qua If $\Gamma$ is a disjoint union of $\Gamma_1$ and a loop, 
then  
  \begin{equation*}
 \bracketG = [3]\bracket{\Gamma_1}
 \end{equation*}

(ii)\qua If  $\Gamma_1$ is obtained by deleting a digon face of 
  $\Gamma$ then 
 \begin{equation*}
  \bracketG = [2]\bracket{\Gamma_1}
 \end{equation*} 

(iii)\qua If  $\Gamma_1, \Gamma_2$ are given by deleting pairs of opposite edges 
 from a square face of $\Gamma$ then 
 \begin{equation*}
 \bracketG= \bracket{\Gamma_1} + \bracket{\Gamma_2}. 
 \end{equation*}

Since a closed web contains either a loop, or a digon face, or a 
square face, these rules determine the invariant. 
For a proof that they are consistent see Kuperberg \cite{Kspiders}. 
 
\begin{figure} [ht!] \drawings{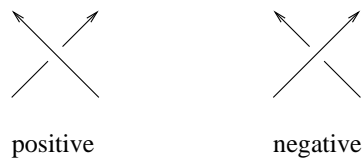}{4.8cm}\caption{Types of crossings} 
 \label{posneg} 
 \end{figure}
Let $D$ be a generic planar projection of an oriented link $L$
and $n$ the number of crossings in $D.$  
We distinguish between \emph{positive} and \emph{negative} crossings, 
see Figure~\ref{posneg}. 

\begin{figure} [ht!] \drawing{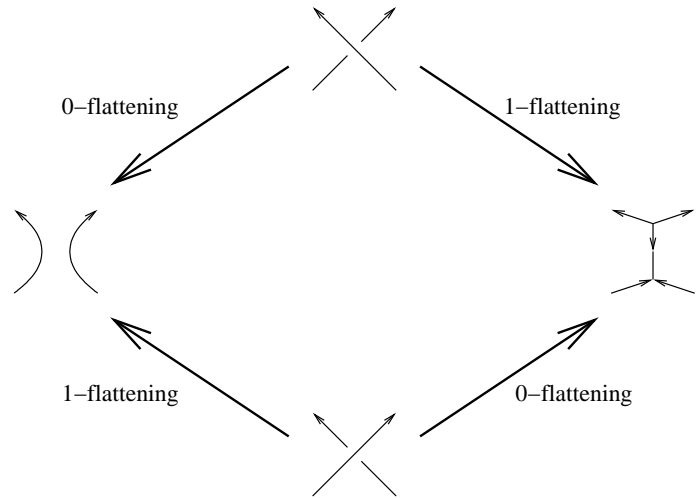}\caption{Flattenings} \label{flats} 
 \end{figure}
\begin{figure} [ht!] \drawings{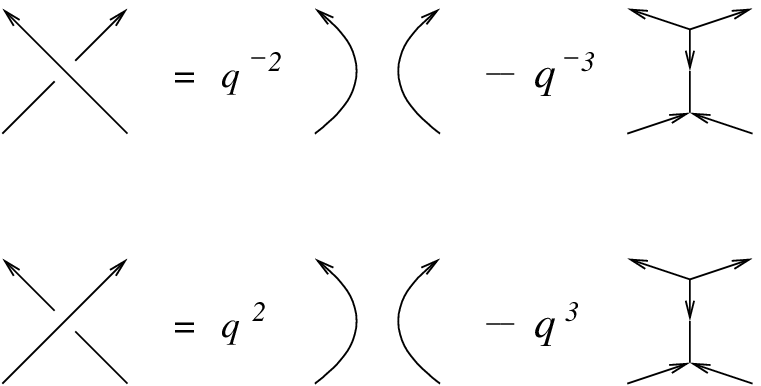}{7.2cm}\caption{Decompositions 
  of crossings} \label{cross} 
 \end{figure}

For each crossing of $D$ consider two  
"flattenings" of this crossing as in Figure~\ref{flats}. 
Define $\bracket{D},$ the bracket of $D,$  as the linear 
combination of the brackets of all $2^n$ flattenings $D_f$ of $D,$ 
via the rules in Figure~\ref{cross},  

\begin{equation*}
\bracket{D} \define \sum_{f\in \{ 0,1\}^n} b_f\bracket{D_f}. 
\end{equation*}
The sum runs over all flattenings $f$ of $D,$ and $b_f$ is 
plus or minus a power of $q,$ see Figure~\ref{cross}. 

\begin{figure} [ht!] \drawings{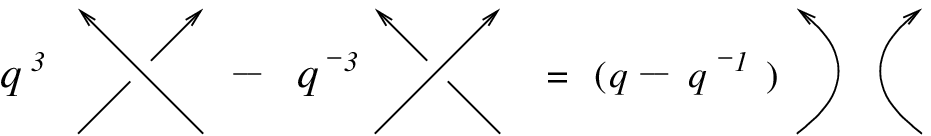}{8.6cm}\caption{Quantum 
$\mf{sl}(3)$ skein formula} \label{homfly3} \end{figure}

\begin{prop} $\bracket{D}$ is preserved by the Reidemeister moves.
\end{prop}
Thus, $\bracket{D}$ is an invariant of an oriented link $L,$ and will 
be denoted $\bracket{L}.$  This invariant of links is associated with the 
quantum group $U_q(\mf{sl}_3)$ and its fundamental representations. 
Excluding rightmost terms from the equations in Figure~\ref{cross} 
yields the $\mf{sl}(3)$ specialization of the HOMFLY skein 
relation (Figure~\ref{homfly3}).  

 $\bracket{L}$
 does not depend on the framing of $L.$ Kuperberg \cite{Kspiders} uses
 a slightly different, frame-dependent, normalization of this
 invariant.


\section{Foams}\label{sfoams}

\subsection{A (1+1)-dimensional TQFT with dots} 
\label{TQFTdots}

Let $\cA= \Z[X] /X^3.$ This ring is commutative Frobenius with 
the trace    
\begin{equation*}\varepsilon: \cA\to \Z, 
 \hspace{0.15in} \varepsilon(X^2) = -1, \hspace{0.15in} \varepsilon(X)=0, 
\hspace{0.15in} \varepsilon(1)=0. 
\end{equation*}
We make $\cA$ into a graded ring, by  
$\deg(1)=-2, \deg(X)= 0, \deg(X^2)=2.$ Multiplication and comultiplication 
in $\cA$ have degree $2,$ the unit and trace maps have degree $-2.$ 
The comultiplication is the map $\cA\to \cA\otimes \cA$ dual to the 
multiplication via the trace form, 
\begin{eqnarray*} 
  \Delta(1) & = & -1\otimes X^2 - X\otimes X - X^2 \otimes 1, \\
  \Delta(X) & = & -X\otimes X^2 - X^2 \otimes X, \\ 
  \Delta(X^2) &  = & -X^2 \otimes X^2.  
\end{eqnarray*}
$\cA$ is the cohomology ring of the complex 
projective plane with the reverse complex orientation, 
\begin{equation*} 
   \cA \cong \mathrm{H}^{\ast}(\overline{\mathbb{P}^2}, \Z), 
\end{equation*}
the trace map being the evaluation on the fundamental $4$-cycle.
The grading of $\mathrm{H}^{\ast}(\overline{\mathbb{P}^2}, \Z)$ 
is shifted downward by $2,$ to balance cohomology groups 
around degree $0.$ Reversing orientation of $\mathbb{P}^2$ 
 changes the trace map, not the multiplication. Note  
 that $(\cA,\varepsilon)$ and $(\cA, - \varepsilon)$ are not 
 isomorphic as Frobenius algebras (but become isomorphic 
 after adding $\sqrt{-1}$ to the ground ring $\Z$). Our 
 choice of trace is essential for what follows. 

The commutative Frobenius ring $\cA$ gives rise to a functor from the 
category of 2-dimensional oriented cobordisms to the category of graded 
abelian groups. We denote this functor by $\cF.$ 
On objects $\cF$ is given by $\cF(M)\cong \cA^{\otimes J}$ 
where $J$ is the set of components of a one-manifold $M.$ To the "pants" 
cobordism from 2 circles to 1 circle $\cF$ associates the multiplication 
map, etc. For more details see \cite[Section 4.3]{BakalovKirillov}, 
\cite{me:jones}. 

All our tensor products are over $\Z$ unless specified otherwise. 
Denote by $\{ n\}$ the shift in the grading up by $n.$

The homomorphism $\cF(S)$ associated with a cobordism $S$ has degree 
$-2 \chi(S)$ where $\chi$ is the Euler characteristic of $S. $
Multiplication by $X$ increases the degree by 2. 

\begin{figure} [ht!] \drawings{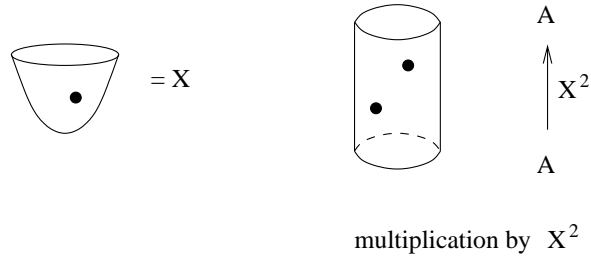}{7.7cm}\caption{The 
 meaning of dots} \label{xandsquare} 
 \end{figure}

A dot on a surface will denote multiplication by $X.$  
For instance, the annulus $S^1\times [0,1]$ 
is the identity cobordism from a circle to itself, and induces the 
identity map $\mathrm{Id}: \cA \to \cA.$ The same annulus with an 
added dot is the multiplication by $X$ endomorphism of $\cA.$ 
Figure~\ref{xandsquare} shows two other examples. The functor $\cF$ 
applied to the "cup" morphism gives the unit map $\Z\to \cA,$ and 
applied to the "cup" with a dot, see Figure~\ref{xandsquare} left,  
produces the map $\Z\to \cA$ which takes $1$ to $X.$ A twice dotted 
annulus is the multiplication by $X^2$ endomorphism of $\cA.$ 
Dots can move freely along a connected component of a surface. 
In the obvious way $\cF$ extends to a functor from the category of 
oriented dotted $(1+1)$-cobordisms to the category of graded abelian groups. 

The degree formula 
 \begin{equation*}
   \mathrm{deg}(\cF(S))= -2\chi(S)
 \end{equation*} 
works for  dotted cobordisms as well if by $\chi(S)$ we denote 
the Euler characteristic of $S$ with dots punctured out. 

The map $\cF(S)$ is zero in each of the following cases: 

(i)\qua a connected component of $S$ has genus $2$ or more, 

(ii)\qua a connected component of $S$ has genus $1$ and at least one 
 dot, 

(iii)\qua a connected component of $S$ contains at least 3 dots.

\begin{figure} [ht!] \drawing{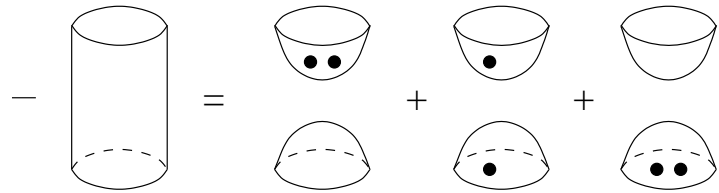}\caption{The surgery formula} 
\label{surgery} 
 \end{figure}

Given any commutative Frobenius algebra $R$ 
with non-degenerate trace $\varepsilon,$ choose a basis $\{ x_i\}$ of $R$ and form 
the dual basis $\{ y_i\}$ determined by $\varepsilon(x_i y_j) = \delta_{i,j}.$ 
Then we have a decomposition $x= \sum y_i \varepsilon (x_i x)$ for any 
$x\in R.$ For $\cA$ and the basis $(1,X,X^2)$ 
the dual basis is $(-X^2, -X,-1).$ The decomposition 
\begin{equation*} 
  -x =  X^2 \varepsilon(x) + X \varepsilon(Xx) + \varepsilon(X^2x)
\end{equation*} 
has a geometric interpretation via a surgery on an annulus, 
see Figure~\ref{surgery}. 

If $S$ is a closed connected oriented surface, possibly decorated 
 by dots, then $\cF(S)=0$ unless 

(i)\qua $S$ is a 2-sphere decorated by 2 dots, then $\cF(M)=-1,$ 

or 

(ii)\qua $M$ is a dotless torus, then $\cF(M)=3.$ 

\begin{figure} [ht!] \drawing{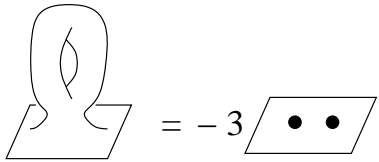}\caption{Genus reduction} 
 \label{torus} 
 \end{figure}

The surgery formula implies the genus reduction formula in Figure~\ref{torus}. 

\subsection{Pre-foams} \label{pre-foams}

\begin{figure} [ht!] \drawing{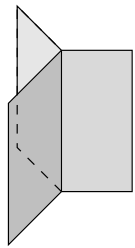}
  \caption{Singularities of a pre-foam} \label{singularities} 
 \end{figure}

{\bf Definition}\qua
Informally, a pre-foam is a 2-dimensional CW-complex with 
singular circles and some additional data.
 A point on a singular circle has a  neighbourhood 
homeomorphic to the product of the letter Y and an interval, 
 see Figure~\ref{singularities}. 

Precisely, a pre-foam $U$ consists of

(i)\qua an orientable smooth compact surface
$U'$ with connected components $U^1, U^2,$  $\dots, U^r.$ 
The boundary components of $U'$ are partitioned into triples 
(the set of boundary components is decomposed into 
a disjoint union of three-element sets). The three 
circles $C^i,C^j,C^k$ in each triple are identified via a 
pair of diffeomorphisms 
  $$ C^i \cong C^j \cong C^k $$
We define the two-dimensional CW-complex underlying $U$ as 
the quotient of $U'$ by these identifications.  

Thus, a triple $(C^i,C^j,C^k)$ of boundary components 
becomes a circle $C$ in the quotient CW-complex $U'/\sim ,$ 
called a \emph{singular} circle. A singular circle 
has a neighbourhood diffeomorphic to 
the direct product of a circle and the letter $Y.$  

We call the images of surfaces $U^i$ in the quotient
CW-complex $U'\hspace{-0.02in}/\hspace{-0.04in}\sim $ 
the \emph{facets} of $U.$   
The image of the interior $U^i\hspace{-0.04in}\setminus 
 \hspace{-0.04in}\partial(U^i)$ of $U^i$  
in $U'\hspace{-0.02in}/\hspace{-0.04in}\sim $ 
 is an open connected orientable surface. 

(ii)\qua A cyclic ordering 
at each singular circle $C$ of the three boundary annuli of $U'$ 
meeting at $C$ (note that there are 
two possible cyclic orderings at each $C$).   

(iii)\qua A number of dots (possibly none) marking each facet 
of $U.$ A dot is allowed to move freely along the facet 
it belongs to, but can't cross a singular circle and move to 
another facet. 

\medskip

{\bf Example}\qua If $U'$ is closed, the pre-foam $U=U'$ is 
a closed surface decorated by dots. 

\medskip
{\bf Example}\qua Let $U'$ be the disjoint union of three discs. 
Then $U$ has one singular circle, and is diffeomorphic to the 2-sphere 
with added equatorial disc. Fix a cyclic ordering of facets 
with the 
equatorial disc followed by the northern hemisphere followed by the 
southern hemisphere followed by the equatorial disc. We call this 
pre-foam the theta-foam. If $a$ dots mark the equatorial disc, 
$b$ dots mark the northern hemisphere, and $c$ dots mark the 
southern hemisphere, we denote this pre-foam by $\theta(a,b,c).$  
The pre-foams $\theta(a,b,c), \theta(b,c,a), \theta(c,a,b)$ 
are isomorphic.  Isomorphism (or homeomorphism) 
of pre-foams is defined in the obvious way.  

\medskip

\sh{Evaluation}
To each pre-foam $U$ we assign an integer $\cF(U),$ the \emph{evaluation} 
of $U,$ via the following rules: 

(1)\qua If $U$ is a closed surface (i.e. has no singular circles), possibly 
with dots, then $\cF(U)$ is defined via the 2D TQFT constructed in  
 section~\ref{TQFTdots}. 

(2)\qua $\cF$ is multiplicative with respect to the disjoint union of pre-foams: 
  \begin{equation*} \cF(U_1\sqcup U_2) = \cF(U_1) \cF(U_2). \end{equation*}
 
(3)\qua Given a circle inside a facet of  
 $U,$ do a surgery on this circle, as in Figure~\ref{surgery}
(on any circle wrapping once around the cylinder on the left 
 hand side). 
 Let $U_1, U_2, U_3$ be the resulting pre-foams (depicted on 
 the right). Then    
 \begin{equation*} 
  - \cF(U) = \cF(U_1) + \cF(U_2) + \cF(U_3).  
 \end{equation*} 
The singular circles of a pre-foam are disjoint from 
the circles on which we can do surgeries. A surgery circle should 
lie in the interior of a facet, and surgeries happen away from singular 
circles. 

(4)\qua The theta-foam evaluates as follows: 
   $$\cF(\theta(a,b,c)) = \begin{cases} 
    0 \mbox{ if }a+b+c\not= 3\mbox{ or }a=b=c=1, \\
    1 \mbox{ if } (a,b,c)=(0,1,2), \mbox{ up to cyclic 
       permutation,} \\
    -1 \mbox{ if } (a,b,c)=(0,2,1), \mbox{ up to cyclic 
       permutation,} \end{cases} $$
 
\begin{figure} [ht!] \drawings{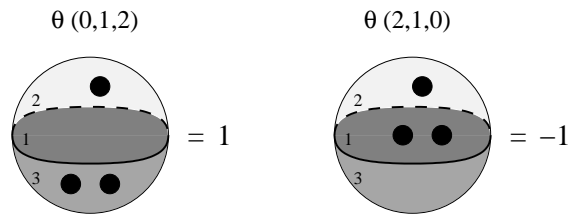}{7.5cm}\caption{Theta-foam 
evaluations. The cyclic order of facets near the 
singular circle is $1\to 2\to 3 \to 1.$}  
  \label{tvalue}
 \end{figure}

\begin{prop} Rules (1)--(4) are consistent and uniquely determine 
$\cF(U)$ for any pre-foam $U.$ 
\end{prop} 

\proof To show that the rules determine the evaluation, 
deform each singular circle $C$ of $U$ to 
three circles $C_1, C_2, C_3$ inside the three annuli of $U$ near $C.$
Do a surgery on each of $C_1, C_2, C_3.$ We get  
 \begin{equation*} 
  \cF(U) = \sum_i a_i \cF(U_i) 
 \end{equation*}
 where $a_i\in \Z,$ and each foam $U_i$ is a disjoint union of closed 
 orientable surfaces and theta-foams, possibly dotted.  Rules (1),(2),(4) 
 determine $\cF(U_i).$ 
 
 To prove consistency, we note that 
 any evaluation via (1)--(4) requires separating each singular 
 circle from the rest of $U$ and creating a theta-foam just for 
 this circle. Functoriality of the functor $\cF$ on dotted surfaces 
(non-singular pre-foams) implies consistency of $\cF(U)$ when 
$U$ has no singular circles. 
Therefore, it suffices to check consistency when 
 $U$ is a theta-foam $\theta(a,b,c).$ We could either evaluate 
 $\theta(a,b,c)$ following (4), or do a surgery on a 
 circle of $\theta(a,b,c)$ that lies in the interior of one 
 of the three discs, and then evaluate using (1),(2), and (4). 
 Direct computations show that these two ways to evaluate 
 a theta-foam give the same answer. 
\endproof

 Rule (4) is related to  the cohomology ring of the full flag 
 variety $Fl_3$ of $\C^3.$ The latter has degree $2$ generators 
 $X_1,X_2,X_3$ and defining relations 
 \begin{eqnarray} 
   X_1 + X_2 + X_3 & = & 0, \label{flag1} \\
   X_1 X_2 + X_1 X_3 + X_2 X_3 & = & 0, \label{flag2} \\
   X_1 X_2 X_3 & = & 0. \label{flag3}  
 \end{eqnarray} 
 The trace form 
 $$\mbox{Tr}: \mathrm{H}^6(Fl_3, \Z) \lra \Z$$ 
 defined by $\mbox{Tr}(X_1 X_2^2) =1$ has 
 the property 
 $$ \mbox{Tr} (X_1^a X_2^b X_3^c) = \cF(\theta(a,b,c)).$$ 
 The dots on the discs of the theta foam correspond to 
  generators $X_1,X_2,X_3$ of the flag variety cohomology 
  ring.  

\begin{prop} If a pre-foam $U_1$ is obtained from a pre-foam $U$ by 
 reversing the cyclic order of the three annuli at a singular circle of $U,$ 
 then 
 \begin{equation*} 
  \cF(U_1) = - \cF(U)
 \end{equation*}
\end{prop}

Define the Euler characteristic $\chi(U)$ of a pre-foam $U$ as the 
Euler characteristic of the two-dimensional CW-complex underlying $U$ 
with all dots punctured out. 

\begin{prop} $\cF(U)=0$ if $\chi(U)\not= 0.$ \end{prop} 

\begin{prop} $\cF(U)=0$ in each of the following cases:   
\begin{itemize} 
\item  a facet of $U$ has genus $2$ or higher, 
\item  a facet of $U$ has genus $1$ and at least one dot, 
\item  a facet of $U$ has at least 3 dots,  
\item  two facets of $U$ that share a singular circle have 
 at least 4 dots,  
\item  two annuli near a singular circle belong to the same 
facet of $U.$ 
\end{itemize} 
\end{prop}

 \sh{Pre-foam identities}

\begin{figure} [ht!] \drawings{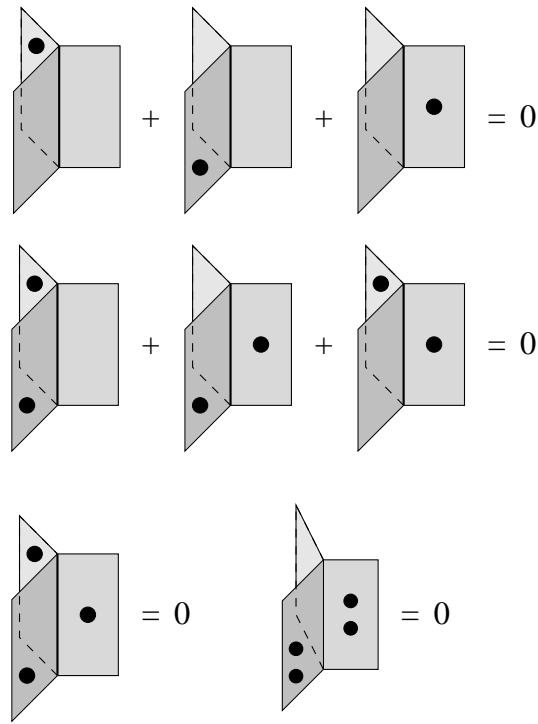}{7cm} \caption{Relations involving 
 dots and singular circles, compare with formulas 
 (\ref{flag1})--(\ref{flag3})}     \label{ssum2}
 \end{figure}

\begin{figure} [ht!] \drawings{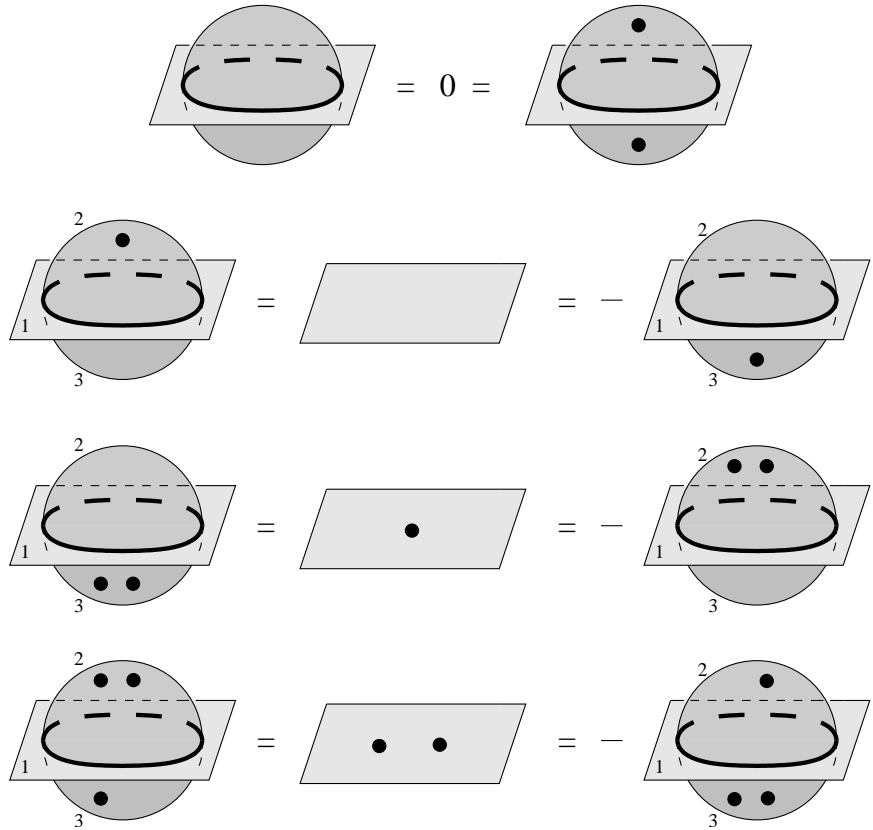}{11.5cm}\caption{Bursting bubbles. 
 Labels $1,2,3$ indicate the cyclic order of facets near a singular 
 circle, $1\to 2 \to 3\to 1$} 
  \label{bubble1} 
 \end{figure}

\begin{figure} [ht!] \drawings{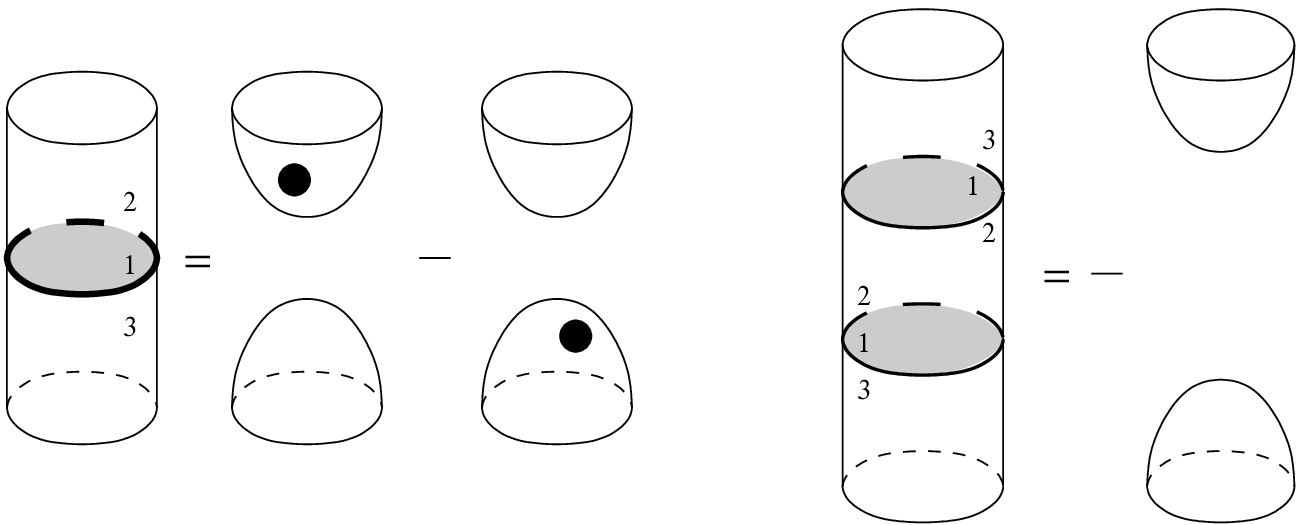}{12.5cm}\caption{Removing discs} 
  \label{tube} 
 \end{figure}

Rules (1)--(4) lead to several handy formulas for pre-foam 
evaluation, depicted in figures~\ref{ssum2}, \ref{bubble1}, 
  \ref{tube}. 

\subsection{Foams and web homology} \label{foams} 

\begin{figure} [ht!] \drawings{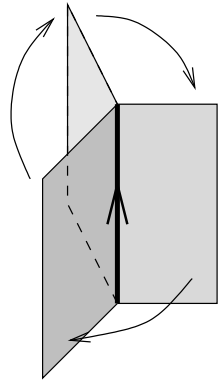}{2.8cm}\caption{Cyclic ordering 
 of annuli determined by singular circle orientation} 
  \label{cyclico} 
 \end{figure}

\sh{Closed foams}
A closed foam is a pre-foam $U$ embedded in $\R^3$ with 
 all facets oriented in such a way that the 
 three annuli near each singular circle are compatibly oriented. 
 Orientations of annuli induce orientations on singular circles. 
 The orientation of a singular circle together with the standard 
 orientation of $\R^3$ induces a cyclic ordering of the three annuli at 
 this circle. We require this ordering to be the same as the 
 ordering specified by the pre-foam $U.$ Figure~\ref{cyclico} explains our 
 cyclic ordering convention.   

There is an obvious forgetful map from closed foams to pre-foams. 
Define the evaluation $\cF(U)\in \Z$ of a closed foam $U$
as the evaluation of the corresponding pre-foam.
  
\medskip

\sh{Foams with boundary} 
Let $T\cong \R^2$ be an oriented plane in $\R^3.$ We say that $T$ intersects 
a closed foam $U$ \emph{generically} if $T\cap U$ is a web in $T.$ Note that the 
orientations of facets of $U$ induce orientations on edges of $T\cap U.$ 

Define a foam with boundary as the intersection of a closed foam  $U$ 
 and $T\times [0,1]\subset \R^3$ such that $T\times \{ 0\}$ 
 and $T\times \{ 1\}$ intersect $U$ generically. We think of 
 a foam $V$ with boundary as a cobordism between webs 
 $\partial_i V \define V \cap T(\times \{ i\})$ for $i=0,1.$ 
 We will always consider the standard embedding 
 of $T\times [0,1]\cong \R^2\times [0,1]$ into $\R^3.$ 
 By a foam we will mean a foam with boundary. 
  We say that two foams are isomorphic if they are isotopic 
 in $\R^2\times [0,1]$ through an isotopy during which the 
 boundary is fixed. 

 \medskip
\textbf{Example}\qua A closed foam is a foam with the empty boundary, 
 and a cobordism between the empty webs. 
\medskip
  
 If $U,V$ are foams and webs $\partial_0 U,\partial_1 V$ 
 are identical as oriented graphs in $\R^2$ (not just isotopic, 
 but isomorphic), we define the composition $UV$ in the obvious 
 way. It is straightforward to  check that any new singular circle 
 that appears in the composition has a neighbourhood 
 homeomorphic to $Y\times \mathbb{S}^1,$ and not to a nontrivial 
 $Y$-bundle over $\mathbb{S}^1.$ Otherwise, the intersection of $UV$ 
with the boundary of a solid torus neighbourhood of that singular circle 
is a simple closed curve $t$ on the 2-torus $T^2$ wrapping three times 
around the longitude and once or twice along the meridian. Thus, the 
curve represents a nontrivial element in $H^1(T^2, \mathbb{Z}_3).$ 
On the other hand, the intersection of $UV$ with the complement 
of the same neighbourhood of the singular circle gives us  
a $\mathbb{Z}_3$ two-chain in the complement with $t$ as its boundary
(to produce the two-chain, contract boundary graphs of $UV$ into 
points). Contradiction. Therefore, $UV$ is a foam. 

 Let $\foams$ be the category with webs as objects and 
 isomorphism classes of foams (rel boundary) as morphisms.    
 A foam $U$ is a morphism from $\partial_0 U $ to $\partial_1 U.$ 
 Define the degree of $U$ as 
 $$ \deg(U) \define \chi(\partial U) - 2\chi(U),$$
   where $\chi$ is the Euler characteristic. 
 $\deg(UV)= \deg(U)+\deg(V)$ for any composable foams $U,V.$

 We now define homology $\cF(\Gamma)$ of a web $\Gamma.$ 

 \medskip

 {\bf Definition}\qua$\cF(\Gamma)$ is an abelian 
 group generated by symbols $\cF(U)$ for all foams $U$ 
 from the empty web $\emptyset$ to 
 $\Gamma.$ A relation $\sum_i a_i \cF(U_i)=0$ for $a_i \in\Z$ and 
 $U_i \in \Hom_{\foams}(\emptyset,\Gamma)$  holds in $\cF(\Gamma)$ 
 if and only if 
 \begin{equation*} 
   \sum_i a_i \cF(V U_i) = 0 
 \end{equation*} 
 for any foam $V$ from $\Gamma$ to the empty web, where 
 $\cF(V U_i)\in \Z$ is the evaluation of the closed foam $VU_i.$  
   
 \medskip

{\bf Example}\qua The homology of the empty web is $\Z.$ 

 \medskip
 
 Since $\deg(VU_i)=\deg(V)+\deg(U_i)$ and  
 $\cF(V U_i)=0$ if $\deg(VU_i)= -2 \chi(VU_i)\not= 0,$ the group  
 $\cF(\Gamma)$ is naturally graded by $\deg(U).$ 

 To a foam $U$ there is associated a homomorphism of abelian 
 groups
 $$ \cF(U): \hspace{0.1in} \cF(\partial_0 U) \lra \cF(\partial_1 U)$$ 
 that takes $\cF(V),$ for any foam 
 $V\in \Hom_{\foams}(\emptyset, \partial_0U),$
 to $\cF(UV).$ This homomorphism has degree $\deg(U).$ 

\begin{figure} [ht!] \drawing{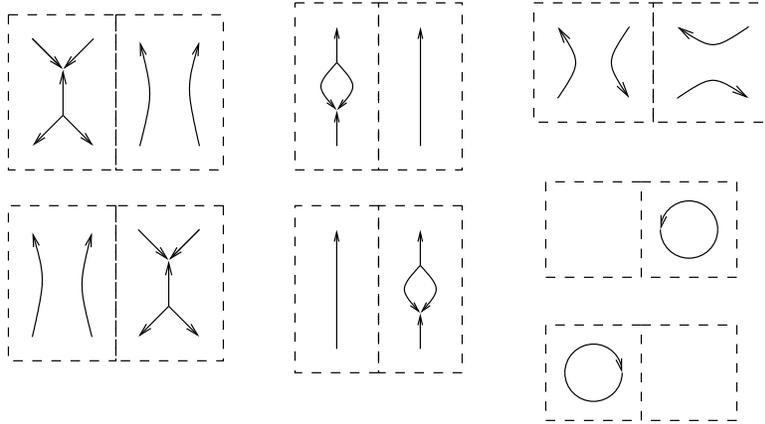}\caption{Basic foams} 
  \label{crossfoam} 
 \end{figure}

We will often depict a dotless foam $U$ by a sequence of its cross-sections 
$U\cap (R^2\times \{ t\})$ for several $t\in [0,1],$ starting 
 with $t=0$ and ending with $t=1.$ Such presentations 
 for foams with only one non-generic cross-section 
 $U\cap (R^2\times \{s\}), s\in [0,1]$ are given in 
 Figure~\ref{crossfoam}, with $t=0,1$ for each foam. We
 call these foams \emph{basic foams}.    
 
\subsection{Web homology skein relations} \label{webskein}

Our goal in this section is to show that $\cF(\Gamma)$ 
 is a free graded abelian group of graded rank $\bracketG$ 
 that satisfies skein relations categorifying those of 
 $\bracketG,$ the latter depicted in Figure~\ref{rules}. 

\medskip 

\sh{Compatibility}
A web $\Gamma$ without trivalent vertices is a closed oriented 
1-manifold embedded in $\R^2.$ Let us show that the two definitions of 
 $\cF(\Gamma)$ for such $\Gamma,$ one in Section~\ref{TQFTdots} and 
 the other in Section~\ref{foams}, give canonically isomorphic abelian 
 groups. To temporarily distinguish them, we call the first  
$\cF_{\ref{TQFTdots}}(\Gamma)$ and the second $\cF_{\ref{foams}}(\Gamma).$ 

\begin{figure} [ht!] \drawings{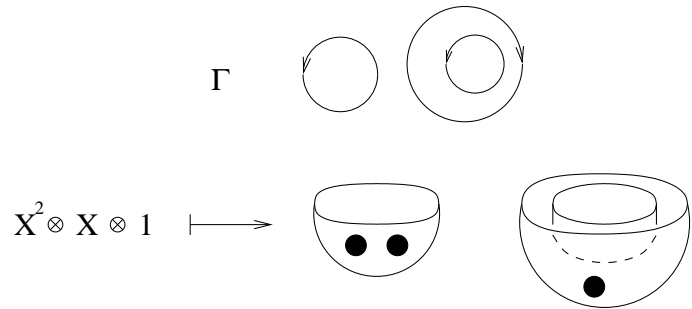}{9cm}\caption{Map $\alpha$} 
  \label{compatible} 
 \end{figure}

Define the following maps: 
$$ \cF_{\ref{TQFTdots}}(\Gamma) \stackrel{\alpha}{\lra}  
   \cF_{\ref{foams}}(\Gamma), \hspace{0.2in} 
  \cF_{\ref{foams}}(\Gamma)\stackrel{\beta}{\lra}
  \cF_{\ref{TQFTdots}}(\Gamma).$$ 
$\cF_{\ref{TQFTdots}}(\Gamma)\cong \cA^{\otimes n},$ 
 where $n$ is the number of connected components of $\Gamma.$
 This abelian group has a basis of elements 
 $X^{a_1}\otimes X^{a_2} \otimes \dots \otimes X^{a_n}$ 
 for  $a_i\in \{ 0,1,2\}.$ We let $\alpha$ take this element to 
 the foam from the empty web to $\Gamma$ which is a union of cups, 
 one for each component of $\Gamma,$ decorated by $a_i$ dots, 
 see example in Figure~\ref{compatible}.  

 To define $\beta,$ start with a foam $U$ from the empty web to $\Gamma,$ 
 and do a surgery near each boundary circle of $U.$ Each term in the 
 resulting sum is a disjoint union of a dotted cup foam and a closed 
 foam. Evaluate closed foams using $\cF$ and convert dotted cup foams 
 to basis elements of $\cF_{\ref{TQFTdots}}(\Gamma)$ by 
 inverting the procedure in Figure~\ref{compatible}. The sum becomes 
 an element of $\cF_{\ref{foams}}(\Gamma),$ and we define 
 $\beta(U)$ as this sum. It is easy to see that $\beta$ is well-defined, 
 and is a two-sided inverse of $\alpha.$

\medskip

\sh{Adding a circle}

\begin{prop} \label{innermost} If $\Gamma$ is obtained from $\Gamma_1$ 
by adding an innermost circle $P$ then  
\begin{equation*} 
  \cF(\Gamma) \cong \cF(\Gamma_1) \otimes \cF(P). 
\end{equation*}  
\end{prop} 
\textbf{Remark}\qua For an example of such $\Gamma$ see Figure~\ref{web}. 
Note that $\cF(P) \cong \cA.$ 

\proof
We leave the proof to the reader. It is similar to the above proof 
that the two definitions of $\cF$ are compatible on 
verticeless $\Gamma.$  \endproof

\begin{figure} [ht!] \drawings{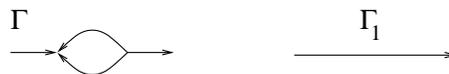}{6cm}\caption{Digon removal} 
  \label{face2a} 
 \end{figure}

\sh{Removing a digon face} 
Suppose that $\Gamma_1$ is given by 
removing a digon face from $\Gamma,$ see Figure~\ref{face2a}. 

\begin{prop} There is a natural isomorphism of graded abelian 
 groups 
\begin{equation} \label{iso2} 
 \cF(\Gamma) \cong \cF(\Gamma_1) \{ 1\} \oplus 
    \cF(\Gamma_1)\{ -1\}. 
 \end{equation} 
\end{prop}  

\begin{figure} [ht!] \drawings{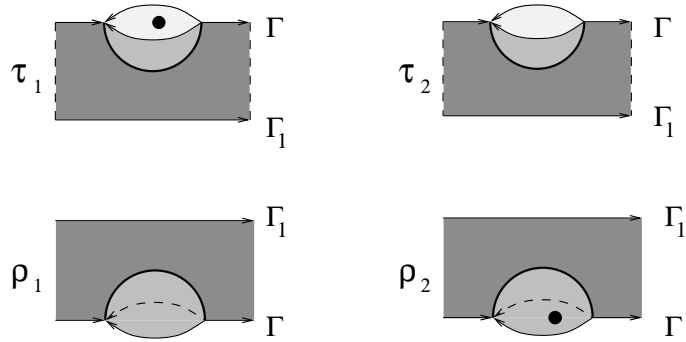}{9cm}\caption{Web cobordisms
 $\tau_1, \tau_2, \rho_1, \rho_2$} 
  \label{face2b} 
 \end{figure}

\proof Consider the cobordisms $\rho_1, \rho_2, \tau_1, \tau_2$ 
between $\Gamma$ and $\Gamma_1$ depicted in Figure~\ref{face2b}.  

These cobordisms induce grading-preserving maps 
 \begin{eqnarray*} 
   \cF(\Gamma_1)\{ 1\} \stackrel{\cF(\tau_1)}{\lra} 
   \cF(\Gamma), & & 
   \cF(\Gamma) \stackrel{\cF(\rho_1)}{\lra} 
   \cF(\Gamma_1)\{ 1\} , \\
   \cF(\Gamma_1)\{ -1\} \stackrel{\cF(\tau_2)}{\lra} 
   \cF(\Gamma),  & & 
   \cF(\Gamma) \stackrel{\cF(\rho_2)}{\lra} 
   \cF(\Gamma_1)\{ -1\}.   
 \end{eqnarray*} 

We claim that 
\begin{eqnarray*} 
    \cF(\rho_1 \tau_1) & = &  \mathrm{Id}_{\Gamma_1} \\
  - \cF(\rho_2 \tau_2) & = &  \mathrm{Id}_{\Gamma_1} \\
    \cF(\rho_1 \tau_2) & = & 0           \\
    \cF(\rho_2 \tau_1) & = & 0           \\
    \cF(\tau_1 \rho_1) - \cF(\tau_2 \rho_2) & = & \mathrm{Id}_{\Gamma}.  
\end{eqnarray*} 
 
\begin{figure} [ht!] \drawings{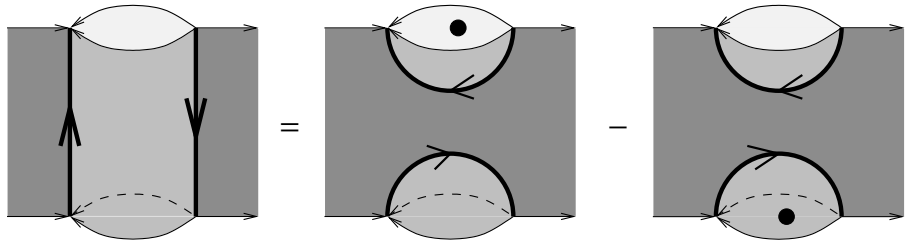}{12cm}\caption{ 
 Identity of $\cF(\Gamma)$ decomposed} 
  \label{face2c} 
 \end{figure}

Indeed, the first four equalities follow from identities for pre-foams, 
 see figures~\ref{ssum2}, \ref{bubble1}. 
The last identity is depicted in Figure~\ref{face2c}. 

To prove it we show that for any foams $U\in \Hom_{\foams}(\emptyset, 
 \Gamma)$ and $V\in \Hom_{\foams}(\Gamma, \emptyset)$ there  is 
 an equality of closed foam evaluations 
 \begin{equation} \label{iden-2-face} 
 \cF(VU)=\cF(V\tau_1 \rho_1 U) - \cF(V \tau_2 \rho_2 U).
 \end{equation}  
There are two cases to consider: 

(1)\qua the two singular intervals of $\mathrm{Id}_{\Gamma}$ shown 
in Figure~\ref{face2c} belong to the same singular circle of 
$VU\cong V\hspace{0.02in}\mbox{Id}_{\Gamma}U.$ 

(2)\qua these intervals belong to different singular circles of $VU.$ 

In each case we do surgeries on $VU$ near all of its singular circles 
that do not intersect the piece of $VU$ shown in Figure~\ref{face2c} left,
and separate these circles from that piece. Identical surgeries are 
 done on the other two closed foams in (\ref{iden-2-face}).
After this operation we are reduced to checking the validity 
 of (\ref{iden-2-face}) in the following cases: 

(1)\qua $VU$ is a theta foam. Foams $V\tau_1 \rho_1 U$ and 
 $V \tau_2 \rho_2 U$ then each have two singular circles. 

(2)\qua $VU$ has two singular circles and is a "connected sum"
 of two theta foams. $V\tau_1 \rho_1 U$ and 
 $V \tau_2 \rho_2 U$ are theta foams.

These cases are verified by hand, using relations in 
figures~\ref{ssum2}, \ref{bubble1}, \ref{tube}.

The five equalities above imply that maps $(\cF(\rho_1), \cF(\rho_2))$ and 
$(\cF(\tau_1), -\cF(\tau_2))$ 
between the left and right hand sides of (\ref{iso2}) 
are mutually inverse isomorphisms. \endproof

\sh{Square face} 

\begin{figure} [ht!] \drawings{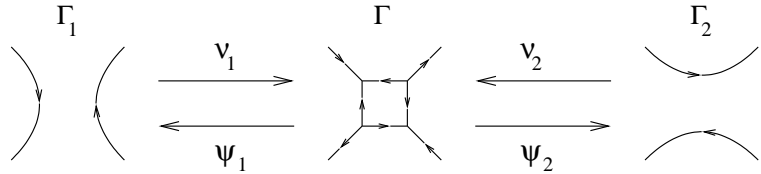}{10cm}\caption{Scheme of cobordisms 
 between $\Gamma$ and its reductions} 
  \label{square1} 
 \end{figure}

\begin{prop} \label{squared} 
 If $\Gamma_1, \Gamma_2$ are two reductions 
 of a web $\Gamma$ with a square face, see Figure~\ref{square1}, 
 there is a canonical isomorphism 
 \begin{equation}\label{sqdecomp}  
   \cF(\Gamma) \cong \cF(\Gamma_1) \oplus \cF(\Gamma_2). 
 \end{equation} 
\end{prop}

\begin{figure} [ht!] \drawings{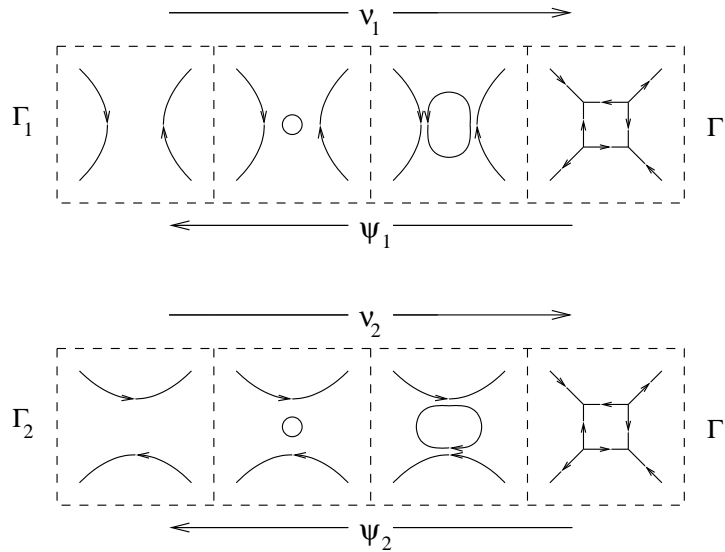}{9.5cm}\caption{
 Cobordisms $\nu_1, \nu_2, \psi_1, \psi_2$} 
  \label{square4} 
 \end{figure}

Let $\nu_1, \nu_2, \psi_1, \psi_2$ be the cobordisms between 
 $\Gamma,\Gamma_1,$ and $\Gamma_2$ depicted in Figure~\ref{square4} 
and denoted by arrows in Figure~\ref{square1}.  

We claim that 
 \begin{eqnarray} 
  \label{eq1} \cF(\psi_1 \nu_1) & = & -\mathrm{Id}_{\Gamma_1} \\
  \label{eq2}  \cF(\psi_2 \nu_2) & = & -\mathrm{Id}_{\Gamma_2} \\
  \label{eq3}  \cF(\psi_2 \nu_1) & = & 0 \\
  \label{eq4}  \cF(\psi_1 \nu_2) & = & 0 \\
  \label{eq5}  
  \cF(\nu_1 \psi_1 ) + \cF(\nu_2 \psi_2) & = & - \mbox{Id}_{\Gamma}. 
  \end{eqnarray}  

Indeed,  (\ref{eq1}) and (\ref{eq2}) are equivalent to 
 the equality in Figure~\ref{tube} right, and 
 (\ref{eq3}) and (\ref{eq4}) are equivalent to  
  the top left equality in Figure~\ref{bubble1}. 

\begin{figure} [ht!] \drawing{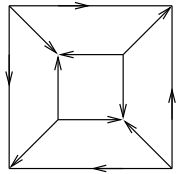}\caption{Cube web $\Sigma$} 
  \label{cubew} 
 \end{figure}

 To prove the last equality, rewrite it as 
 $$\cF(\nu_1 \psi_1 ) + \cF(\nu_2 \psi_2) + \mbox{Id}_{\Gamma}=0.$$ 
 We need to show that, for any $V\in \Hom_{\foams}(\emptyset, \Gamma)$ 
 and $U\in \Hom_{\foams}(\Gamma, \emptyset),$ the following relation 
 between three integers holds:
 $$\cF(U\nu_1 \psi_1V) + \cF(U\nu_2 \psi_2 V) + \cF(UV)=0.$$ 
 Non-trivial parts of web cobordisms $\nu_1\psi_1$ and $\nu_2 \psi_2$ 
 lie in $D^2\times [0,1],$ the direct product of a 2-disc and the 
 interval $[0,1].$ Smooth out the sharp edges of $D^2\times [0,1]$ 
 slightly so that its boundary is a 2-sphere in $\R^3.$ Denote this 
 sphere by $S$ and the 3-disc that it bounds by $D^3.$ 
 Each of the cobordisms $\nu_1\psi_1, \nu_2 \psi_2,$ and 
 $\mbox{Id}_{\Gamma}$ intersect $S$ along the same web $\Sigma,$
 isomorphic to the cube web, Figure~\ref{cubew}. 

\begin{figure} [ht!] \drawings{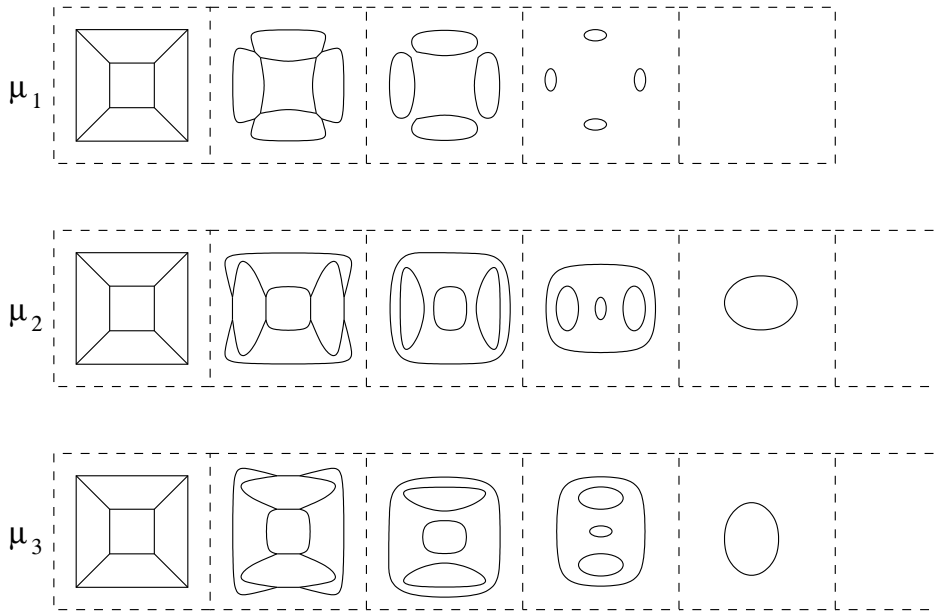}{12.5cm}\caption{
 Cobordisms $\mu_1, \mu_2, \mu_3$ from $\Sigma$ to the 
 empty web. Orientations of edges are omitted to avoid clutter.} 
  \label{square6} 
 \end{figure}

 Each of the three cobordisms $U\nu_1\psi_1 V, U\nu_2 \psi_2V $ and 
 $UV\cong U\hspace{0.02in}\mbox{Id}_{\Gamma}V$ has 
 identical intersection, 
 denoted $E,$  with the complement of $D^3$ in $\R^3.$
 This intersection can be viewed as a foam whose boundary is 
 $\Sigma.$ Apply a diffeomorphism of $\mathbb{S}^3= \R^3
 \cup \{ \infty\}$ that takes $S$ to $\R^2\cup \{ \infty\}.$ 
 This diffeomorphism takes $\Sigma$ to a planar web, also denoted 
 $\Sigma,$ and $E$ to a foam in $\R^3$ with boundary $\Sigma.$ 
 The images $\mu_1, \mu_2, \mu_3$ of intersections 
 of $U\nu_1\psi_1 V, U\nu_2 \psi_2V $ and $UV$ 
 with $D^3$ are foams with boundary $\Sigma,$ 
 depicted in Figure~\ref{square6}. 

 This transformation preserves the values of all closed foams. 
 Therefore, it suffices to show that 
 \begin{equation} \label{square-close} 
 \cF(\mu_1 E) + \cF(\mu_2 E) + \cF(\mu_3 E)=0
 \end{equation} 
 for any foam $E\in \Hom_{\foams}(\emptyset, \Sigma).$ 
 The foam $E$ contains four singular arcs that pair up the eight 
 vertices of $\Sigma,$ each arc connecting an "in" vertex with 
 an "out" vertex. 

 Suppose that one of these arcs connects two neighbouring 
 vertices $v_1, v_2$ 
 of $\Sigma.$ One of the three foams $\mu_1, \mu_2, \mu_3$ 
 contains a singular arc connecting the same pair of vertices. 
 Since $\Sigma$ has automorphisms that extend to cyclic 
 permutations of $\mu_1, \mu_2, \mu_3,$ we can assume that 
 this foam is $\mu_1.$ Simplify closed foams 
 $\mu_1 E, \mu_2 E, \mu_3 E$ by cutting off singular circles 
 that are disjoint from $\Sigma.$ The foam $\mu_1 E$ has a facet 
 containing the edge of $\Sigma$ that joins vertices $v_1$ and $v_2.$ 
 Using relations in figures~\ref{torus} and \ref{ssum2},  
 reduce the genus of this facet to $0$ and move all dots away 
 from it to other parts of $E.$ Now apply the relation in Figure~\ref{tube} 
 left to the resulting foam, and the Figure~\ref{face2c} equality 
 at certain digon sections of $\mu_2E$ and $\mu_3E.$ Everything in 
 the left hand side of (\ref{square-close}) cancels and identity 
 (\ref{square-close}) follows in this case. 

 The only other case to consider is when each of the four arcs of $E$ 
 connect a vertex of $\Sigma$ with the opposite vertex (the one 
 that is three edges away). To treat this case, do 
 surgery on $E$ to remove singular circles that do not 
 intersect $\Sigma$ and reduce the genus of all facets of $E$ to $0.$ 
 The Euler characteristic of $\mu_iE$ is then $-4 + 2m$ where $m$ is 
 the number of dots on (modified) $E.$ Unless $E$ has exactly two 
 dots, pre-foams $\mu_1 E, \mu_2 E, \mu_3E$  
 have non-zero Euler characteristic and their evaluations equal $0.$ 
 Furthermore, ignoring the dots, pre-foams $\mu_1 E, \mu_2 E, 
 \mu_3 E$ are isomorphic (these isomorphisms come from 
 automorphisms of $\Sigma$ that extend to modified $E,$ the 
 latter viewed as a pre-foam with boundary). Move the two 
 dots close to each other and to a vertex of $\Sigma.$ 
 Relation (\ref{square-close}) now follows from the top two 
 relations in  Figure~\ref{ssum2}.  
 
 We have shown that 
 $$\cF(\mu_1) + \cF(\mu_2) + \cF(\mu_3) = 0,$$ 
 which is equivalent to (\ref{eq5}).

 Equations (\ref{eq1}) through (\ref{eq5}) imply that $(\cF(\psi_1), \cF(\psi_2))$ and 
  $(-\cF(\nu_1), -\cF(\nu_2))$ are mutually inverse isomorphisms 
 between the L.H.S. and the R.H.S. of (\ref{sqdecomp}). 
 Proposition~\ref{squared} follows. \endproof

\begin{corollary} \label{cor-rank}
$\cF(\Gamma)$ is a free abelian group of graded 
rank $\bracketG.$ 
\end{corollary}

\begin{corollary} There is a natural isomorphism   
 \begin{equation*} 
  \cF(\Gamma_1 \sqcup \Gamma_2 ) \cong \cF(\Gamma_1) \otimes 
  \cF(\Gamma_2), 
 \end{equation*}
where $\Gamma_1 \sqcup \Gamma_2$ is the disjoint union of 
 webs $\Gamma_1, \Gamma_2.$  
\end{corollary}


\section{Chain complex of a link diagram} \label{chains} 

Given a plane diagram $D$ of a link $L,$ let $I$ be the set of 
crossings of $D,$ and $p_+,$ respectively $p_-,$  
be the number of positive, respectively negative, crossings;
 $p_++p_-=|I|.$ 

To $D$ we associate an $|I|$-dimensional cube of graded 
abelian groups. Vertices of the cube are in a bijection with 
subsets of $I.$ To $J\subset I$ we associate a web $D_J,$ which is 
the $J$-flattening of $D.$ If $a\in J,$ the crossing $a$ gets 
$1$-flattening, otherwise it gets $0$-flattening, see Figure~\ref{flats}.
 In the vertex $J$ of the cube we place the graded 
abelian group $\cF(D_J)\{3p_- - 2p_+ -|J|\},$ and
 to an inclusion $J\subset Jb,$ where $b\in I\setminus J,$ 
and $Jb$ denotes the set $J\sqcup \{ b\},$ we assign 
the homomorphism 
$$ \cF(D_J)\{3p_- - 2p_+ -|J|\} \lra 
 \cF(D_{Jb})\{3p_- - 2p_+ -|J|-1\}$$ 
induced by the basic cobordism between webs $D_J$ and $D_{Jb},$
which looks like the cobordism in Figure~\ref{basic} 
 (if crossing $b$ is negative) or its reflection in a horizontal 
 plane (if crossing $b$ is positive).  
Each of these homomorphisms is grading-preserving. 
Functoriality of $\cF$ implies that the cube is commutative. 

We add minus signs to some maps to make each square in the cube 
anticommute. An intrinsic way to assign minus signs can be 
found, for example, in \cite[Section 3.3]{me:jones}.
Next form the total complex of this anticommutative cube, 
placing the first non-zero term, $\cF(D_{\emptyset})\{3p_--2p_+\},$ 
in cohomological degree $-p_-.$ Denote the total complex 
by $\cF(D).$ It is non-zero in cohomological degrees between 
$-p_-$ and $p_+.$ The term of degree $-p_-+i$ is the following 
direct sum: 
$$\cF(D)^{-p_-+i} = \oplusop{J\subset I, |J|=i} \cF(D_J) 
 \{3p_- - 2p_+ -i \}.$$ 
The differential in the complex preserves the grading. 

\begin{figure} [ht!] \drawings{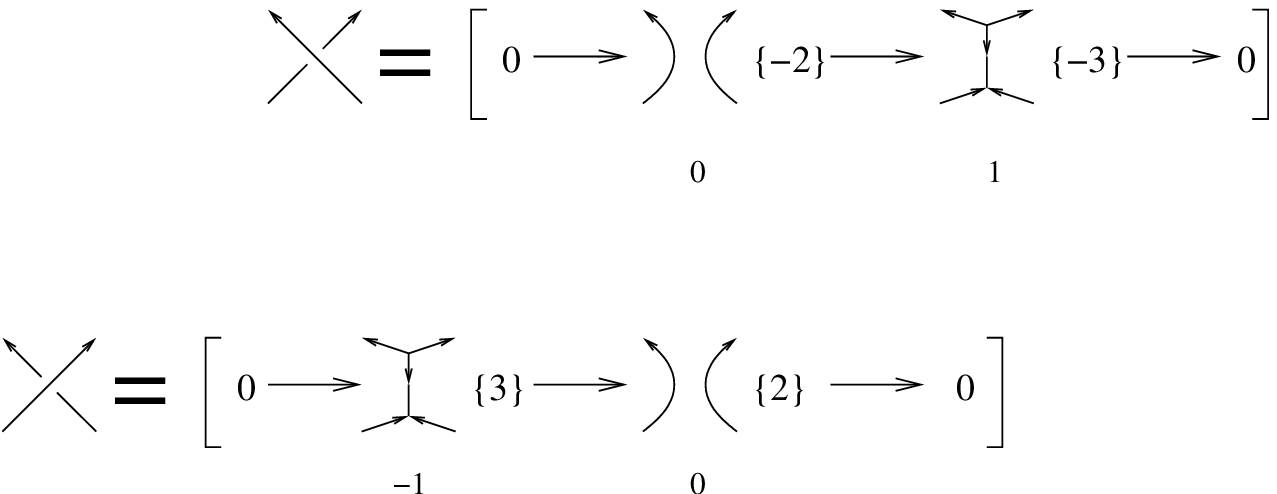}{11.5cm}\caption{$-1,0,1$ 
 underlining the flattenings indicate cohomological degree} \label{cones} 
 \end{figure}
Figure~\ref{cones} explains our construction of $\cF(D).$ 
Corollary~\ref{cor-rank} implies that the Euler characteristic 
of $\cF(D)$ is the Kuperberg bracket $\bracket{L}.$ In the next section 
we show that $\cF(D),$ up to chain homotopy equivalence, is 
 a link invariant.


\section{Invariance under Reidemeister moves} \label{invariance}

Proofs of all propositions stated in this section are straightforward 
 and left to the reader. 

\subsection{Type I moves} 

\begin{figure}[ht!]  \drawings{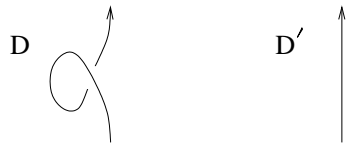}{4.5cm}\caption{A type I move} 
  \label{curla1} 
 \end{figure}

We kick-off by showing that complexes $\cF(D)$ and $\cF(D')$ are 
homotopy equivalent, for diagrams $D$ and $D'$ depicted in 
Figure~\ref{curla1}. 

\begin{figure}[ht!]  \drawings{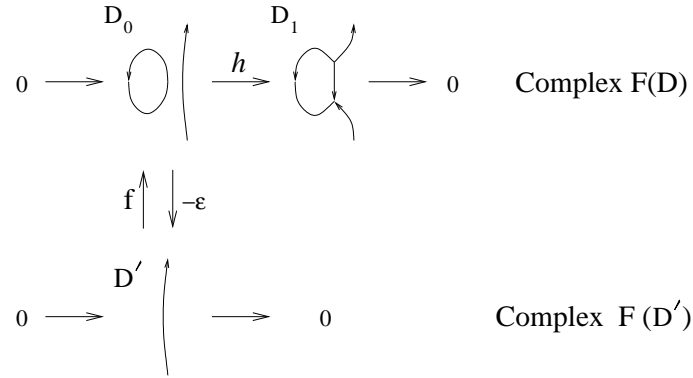}{9cm}\caption{Realization of 
 $\cF(D')$ as a direct summand of $\cF(D)$} 
  \label{curla2} 
 \end{figure}

Denote the flattenings of $D$ by $D_0$ and $D_1$. 
The complex $\cF(D)$ is the total complex of the bicomplex  
$$0 \lra \cF(D_0)\{-2\} \stackrel{h}{\lra} \cF(D_1)\{ -3\}\lra 0, $$ 
see the top of Figure~\ref{curla2}. The map $h$ is induced by 
the basic cobordism from $D_0$ to $D_1.$ 

Define $f: \cF(D') \lra \cA\otimes \cF(D') \cong \cF(D_0)$ 
by 
$$f(a) = X^2 \otimes a + X \otimes Xa + 1 \otimes X^2 a. $$  
Let $\varepsilon: \cF(D_0) \to \cF(D')$ be the trace map, 
 $$\varepsilon(x\otimes a) = \varepsilon(x) a,\hspace{0.1in} 
 \mathrm{for} \hspace{0.1in}x\otimes a\in \cA\otimes \cF(D')\cong \cF(D_0).$$
From the formulas $ -\varepsilon f = \mbox{Id}$ and 
 $\mbox{im}(f) = \mbox{ker}(h)$ 
we get a direct sum decomposition 
$$\cF(D)\cong f(\cF(D')) \oplus \{ 0 \to \mbox{ker}(\varepsilon) 
 \to \cF(D_1) \to 0 \}. $$ 
The complex $f(\cF(D'))$  is isomorphic to $\cF(D'),$ the second 
 direct summand is contractible, since $h,$ restricted to the kernel 
of $\varepsilon,$ is an isomorphism of complexes. Therefore, $\cF(D)$ and 
 $\cF(D')$ are homotopy equivalent.

\begin{figure}[ht!]  \drawings{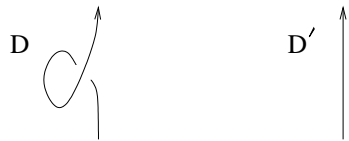}{4.5cm}\caption{The other type I move} 
  \label{curlb1} 
 \end{figure}

We next take care of the other Reidemeister type I move, see 
 Figure~\ref{curlb1}.

\begin{figure}[ht!] \drawings{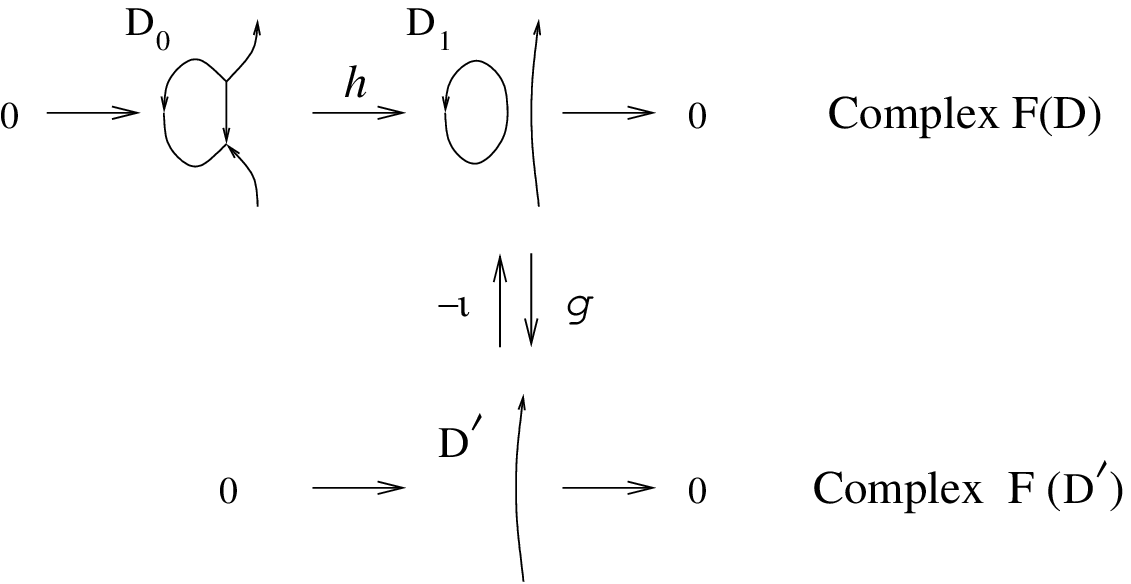}{9cm}\caption{$\cF(D')$ as a 
 direct summand of $\cF(D)$} 
  \label{curlb2} 
 \end{figure}

 The complex $\cF(D)$ is the total complex of the bicomplex 
$$0 \lra \cF(D_0)\{3\} \stackrel{h}{\lra} \cF(D_1)\{ 2\}\lra 0, $$ 
see the top of Figure~\ref{curlb2}. 

Define $g: \cF(D_1)\cong \cA \otimes \cF(D')\lra \cF(D')$ by
$$g(x\otimes a) = \varepsilon(x) X^2a + \varepsilon(Xx) Xa + 
 \varepsilon(X^2x)a, \hspace{0.1in} x\otimes a\in \cA\otimes \cF(D'),$$
and $\iota: \cF(D') \lra \cA \otimes \cF(D') \cong \cF(D_1)$ 
is the unit map, $\iota(a) = 1 \otimes a.$  

Since $-g \iota=\mbox{Id}, g h=0$ 
 and $\cF(D_1) \cong h(\cF(D_0)) \oplus 
 \iota(\cF(D')),$ there is a direct sum decomposition 
$$\cF(D)\cong \iota(\cF(D')) \oplus \{ 0 \to \cF(D_0) \stackrel{h}{\to} 
 \mbox{ker}(g) \to 0 \}. $$
Complex $\iota(\cF(D'))$ is isomorphic to $\cF(D'),$ while 
the second direct summand is contractible,  
since $h$ is an isomorphism onto $\mbox{ker}(g).$  
This establishes a homotopy equivalence between $\cF(D)$ and $\cF(D').$ 

\subsection{Type II moves} 

\begin{figure}[ht!]  \drawings{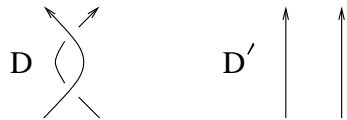}{4.5cm}\caption{A type II move} 
  \label{ta1} 
 \end{figure}

\begin{figure}[ht!]  \drawings{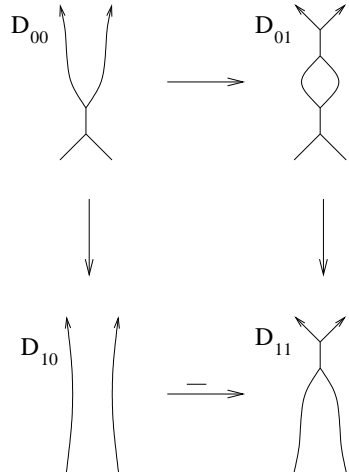}{4.5cm}\caption{Flattenings of D} 
  \label{ta2} 
 \end{figure}

Consider diagrams $D$ and $D'$ in Figure~\ref{ta1}. Figure~\ref{ta2} 
depicts the four flattenings of $D.$
Arrows correspond to basic cobordisms, and 
a bar over the bottom arrow indicates that the map induced by 
this cobordism will be taken with the minus sign when the total 
complex of the commutative square in Figure~\ref{ta2} is formed.

\begin{figure} [ht!] \drawings{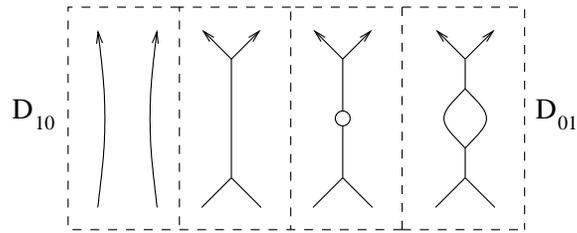}{7.5cm}\caption{Cobordism $\alpha$} 
  \label{ta3} 
 \end{figure}

\begin{figure} [ht!] \drawings{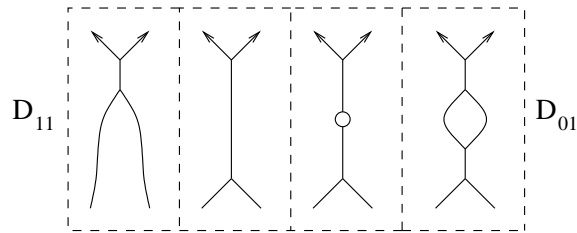}{7.5cm}\caption{Cobordism $\beta$} 
  \label{ta4} 
 \end{figure}

Let $\alpha$ be the cobordism depicted in Figure~\ref{ta3} and 
 $\beta$ the cobordism in Figure~\ref{ta4}. 

Let 
$$Y_1=\{(x, \cF(\alpha)x) \subset \cF(D_{10}) \oplus \cF(D_{01}),
 \hspace{0.1in} x\in \cF(D_{10})\}.$$
Let $Y_2$ be the subcomplex of $\cF(D)$ generated by $\cF(D_{00}).$ 
Let 
$$Y_3= \{(\cF(\beta)x, y) \subset \cF(D_{01})\oplus \cF(D_{11}), 
  \hspace{0.1in} x,y\in \cF(D_{11})\}.$$ 
$Y_3$ is the subcomplex of $\cF(D)$ generated by $\cF(\beta)x, $ 
over all $x\in \cF(D_{11}).$

\begin{prop} \begin{enumerate} \item 
 The complexes $Y_2$ and $Y_3$ are contractible. 
 \item $Y_1$ is a subcomplex of $\cF(D)$ isomorphic to $\cF(D').$ 
 \item $\cF(D) \cong Y_1 \oplus Y_2 \oplus Y_3.$
  \end{enumerate} 
\end{prop} 
Therefore, $\cF(D)$ and $\cF(D')$ are homotopy equivalent.

\begin{figure} [ht!] \drawings{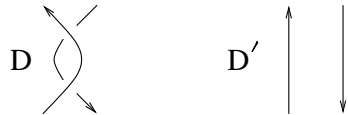}{4.5cm}\caption{Another type II move} 
  \label{tb1} 
 \end{figure}

We next consider a type II move with a different orientation, 
see Figure~\ref{tb1}. 

\begin{figure} [ht!] \drawings{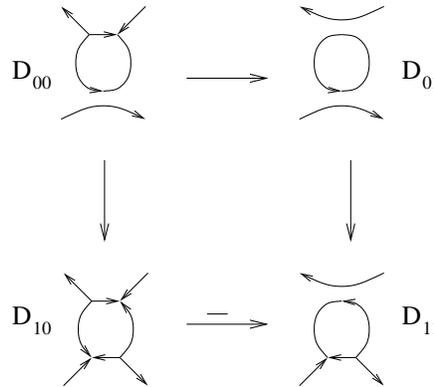}{5.7cm}\caption{Flattenings of $D$} 
  \label{tb2} 
 \end{figure}

The four flattenings of $D$ are depicted in Figure~\ref{tb2}. 

\begin{figure} [ht!] \drawings{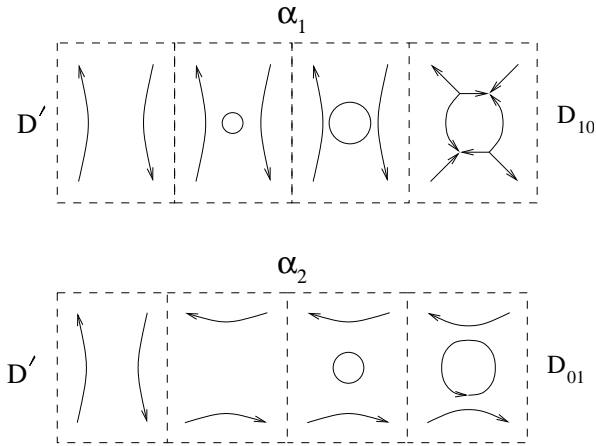}{7.8cm}\caption{Cobordisms $\alpha_1$ 
 and $\alpha_2$} 
  \label{tb3} 
 \end{figure}

Let 
$$Y_1= (\cF(\alpha_1) x, \cF(\alpha_2) x) \subset \cF(D_{10})
   \oplus \cF(D_{01}), $$
 where $x\in \cF(D')$ and $\alpha_1, \alpha_2$ are Figure~\ref{tb3}
 cobordisms.  

Let $Y_2$ be the subcomplex of $\cF(D)$ generated by $\cF(D_{00}).$

The diagram $D_{01}$ contains  a circle, see Figure~\ref{tb2} top 
 right. Let $D''$ be the diagram given by removing this circle 
 from $D_{01}.$ Then  $\cF(D_{01})\cong \cA\otimes \cF(D'').$ 
Let $Y_3$ be the subcomplex of $\cF(D)$ generated by
$1\otimes \cF(D'')\subset \cF(D_{01})\subset \cF(D)$ and 
$X\otimes \cF(D'')\subset \cF(D_{01})\subset \cF(D).$ 

\begin{prop} \begin{enumerate} 
\item Complexes $Y_2$ and $Y_3$ are contractible. 
\item $Y_1$ is a subcomplex of $\cF(D)$ isomorphic to $\cF(D').$ 
\item $\cF(D) \cong Y_1 \oplus Y_2 \oplus Y_3.$ 
\end{enumerate} 
\end{prop} 

The proposition implies that $\cF(D)$ and $\cF(D')$ are 
homotopy equivalent. 

\subsection{Type III move} 

\begin{figure} [ht!] \drawings{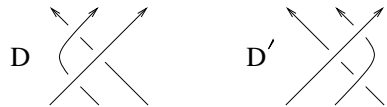}{5cm}\caption{Reidemeister III move} 
  \label{braid1} 
 \end{figure}

Any two moves of type III are equivalent modulo type I and II moves. 
Therefore, it suffices to treat only one case of the type III moves. 
We choose the one in Figure~\ref{braid1}. 

\begin{figure} [ht!] \drawing{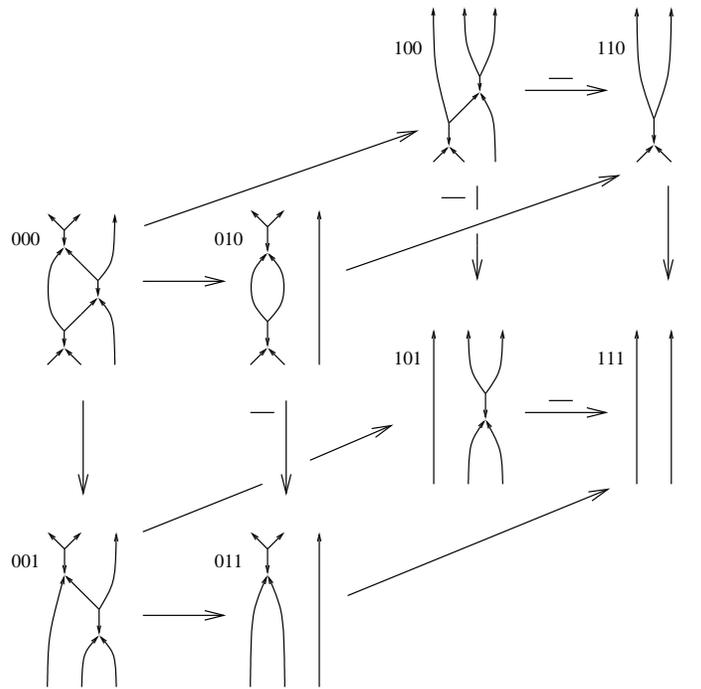}\caption{The cube of eight 
 flattenings of $D$} 
  \label{braid2} 
 \end{figure}

$\cF(D)$ is the total complex of the cube of $8$ complexes shown in 
Figure~\ref{braid2}.  

\begin{figure} [ht!] \drawing{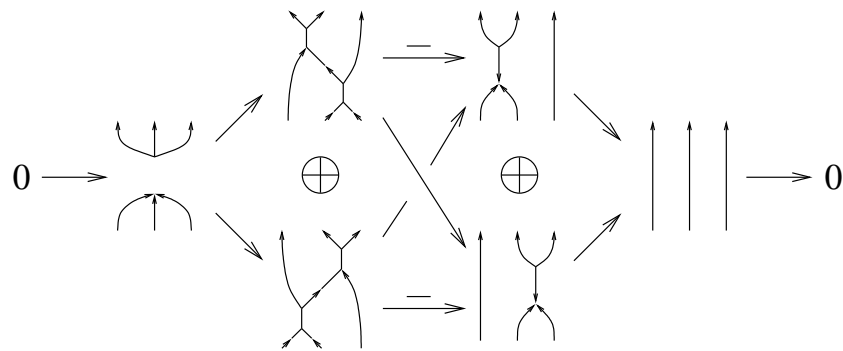}\caption{Complex $Q$} 
  \label{compact} 
 \end{figure}

Let $Q$ be the complex associated to Figure~\ref{compact}. It is 
the total complex 
of the bicomplex built out of complexes assigned to the six diagrams
 in Figure~\ref{compact}.  Five of these diagrams 
are $D_{001},$  $D_{100},$  $D_{101},$  $D_{110},$ $D_{111}.$ The maps 
are induced by basic cobordisms between the diagrams. 

\begin{figure}[ht!] \drawings{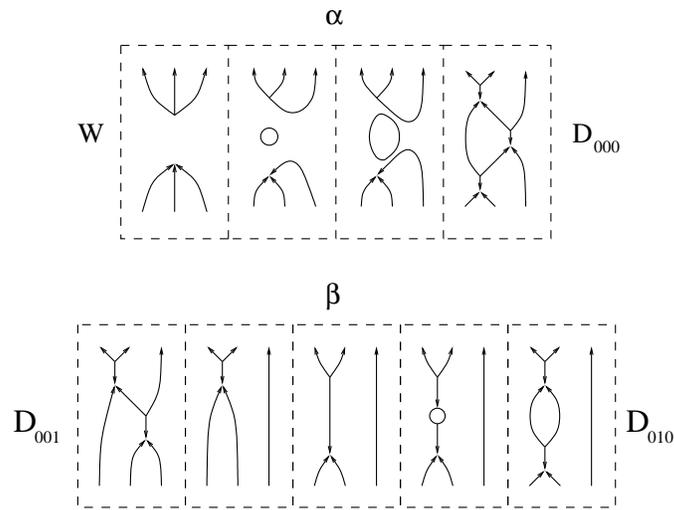}{8.8cm}\caption{Cobordisms $\alpha$ 
   and $\beta$} \label{braid4} 
 \end{figure}

Let $Y_1$ be the subcomplex of $\cF(D)$ which, as an abelian group, 
is the direct sum of 
\begin{itemize} 
\item $\cF(D_{1ij}),$ over $i,j\in \{ 0,1\},$ 
\item $\cF(W)\subset \cF(D_{000})$ where $W$ is the leftmost 
 diagram in Figure~\ref{compact}, and the inclusion is induced 
 by cobordism $\alpha$ in Figure~\ref{braid4}. 
\item $(x, \cF(\beta) x),$ for $x\in \cF(D_{001})$ and cobordism 
 $\beta$ depicted in Figure~\ref{braid4}.  
\end{itemize} 

\begin{figure} [ht!] \drawings{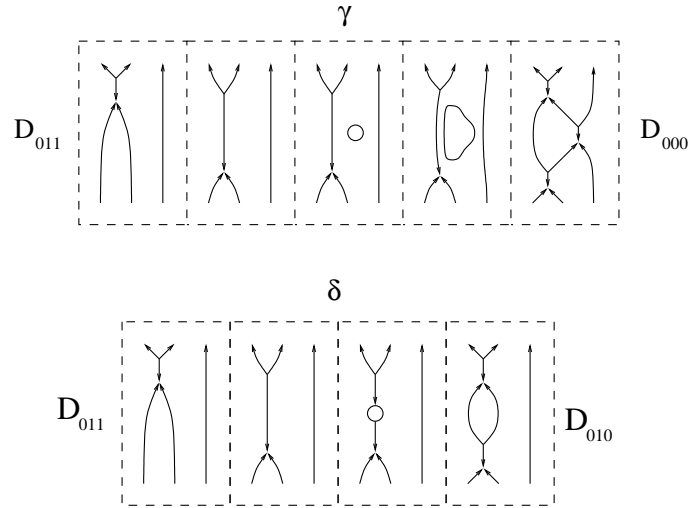}{9cm}\caption{Cobordisms $\gamma$ 
   and $\delta$} \label{braid5} 
 \end{figure}

Let $Y_2$ be the subcomplex of $\cF(D)$ generated by  
the subcomplex of $\cF(D_{000})$ isomorphic to $\cF(D_{011})$ 
via the inclusion induced by cobordism $\gamma$ in Figure~\ref{braid5}. 

Let $Y_3$ be the subcomplex of $\cF(D)$ generated by 
the subcomplex of $\cF(D_{010})$ isomorphic to $\cF(D_{011})$ 
via the inclusion induced by cobordism $\delta$ in Figure~\ref{braid5}. 
 
\begin{prop} $\cF(D)\cong Y_1\oplus Y_2 \oplus Y_3.$ The complexes 
 $Y_2$ and $Y_3$ are contractible. The complexes $Y_1$ and $Q$ are 
 isomorphic. 
\end{prop} 
A similar decomposition establishes homotopy equivalence of 
 $\cF(D')$ and $Q,$ and implies that $\cF(D)$ and $\cF(D')$ are 
 homotopy equivalent complexes.

\section{Flag varieties} \label{flagvar} 

To a closed web $\Gamma$ associate a commutative graded ring $R(\Gamma)$ with 
generators $X_i,$ over all edges $i$ of $\Gamma,$ and relations 
\begin{equation}\label{symmetric}
 X_i + X_j + X_k =0, \hspace{0.1in} X_i X_j + X_i X_k + X_j X_k =0, 
\hspace{0.1in} X_i X_j X_k =0. 
\end{equation}
 whenever edges $i,j,k$ meet at a vertex, and 
$X_i^3=0$ if $i$ is a closed loop. We set the degree of each $X_i$ to $2.$ 
Relations (\ref{symmetric}) are equivalent to the condition that every symmetric 
polynomial in $X_i,X_j,X_k$ vanishes. Relation $X_i^3=0$ for any 
non-closed edge $i$ follows from (\ref{symmetric}).  

$\cF(\Gamma)$ is a module over $R(\Gamma).$ The cobordism 
of merging a circle near edge $i$ with this edge induces a map 
 $\cA \otimes_{\Z} \cF(\Gamma) \lra \cF(\Gamma).$ Multiplication by 
$X\in \cA$ is then an endomorphism $X_i$ of $\cF(\Gamma)$ of degree $2.$ 
These endomorphisms satisfy relations (\ref{symmetric}) and make 
$\cF(\Gamma)$ a module over $R(\Gamma).$ 

If $\Gamma$ can be reduced to a circle by removing digon faces, 
see example in Figure~\ref{theta3} left  
(we call such $\Gamma$ a \emph{digon} web),  
then $\cF(\Gamma)$ is a free 
module over $R(\Gamma)$ of rank one. If $\Gamma$ is the cube web, 
 see Figure~\ref{theta3} right, then $\cF(\Gamma)$ 
is not free over $R(\Gamma).$

Let $V$ be a three-dimensional complex vector space equipped with a 
hermitian form. Let $P(V)$ be the projective space of complex lines in $V.$ 
To $V$ and $\Gamma$ we associate a topological space $V(\Gamma),$ a 
subspace of $P(V)^{\times E},$ where $E$ is the set of edges of $\Gamma.$ 
A set of lines $\{ L_a\}_{a \in E},$  $L_a\subset V,$ is in $V(\Gamma)$ 
iff $L_a$ and $L_b$ are orthogonal whenever $a$ and $b$ share a common vertex. 

The inclusion $V(\Gamma) \subset P(V)^{\times E} $ induces a homomorphism 
of cohomology rings $H(P(V)^{\times E},\Z) \to H(V(\Gamma), \Z)$ which factors 
through to a homomorphism $u: R(\Gamma) \to H(V(\Gamma), \Z).$ If 
$\Gamma$ is a digon web, $u$ is a ring isomorphism. We don't know if 
$u$ is an isomorphism for any $\Gamma.$ 

If $\Gamma$ is a digon web, 
 $$ \cF(\Gamma) \cong \mathrm{H}^{\ast}(V(\Gamma), \Z)\{-\frac{n}{2}-2\}.$$


\Addresses\recd

\end{document}

%% file: agtout.tex
\def\ifplaintex{\expandafter\ifx\csname documentclass\endcsname\relax}

\def\gtp{{\mathsurround=0pt\it $\cal G\mskip-2mu$eometry \&\ 
$\cal T\!\!$opology $\cal P\!$ublications}}  

\def\recd{{\small Received:\qua\receiveddate\ifx\reviseddate\relax
\else\qquad Revised:\qua\reviseddate\fi\par}} 

\def\lognumber#1{\def\thelognumber{#1}}
\def\volumenumber#1{\def\thevolumenumber{#1}}
\def\volumeyear#1{\def\thevolumeyear{#1}}
\def\papernumber#1{\def\thepapernumber{#1}}
\def\pagenumbers#1#2{\def\startpage{#1}\def\finishpage{#2}}
\def\published#1{\def\publishdate{#1}}

\def\received#1{\def\receiveddate{#1}}
\def\revised#1{\def\reviseddate{#1}}
\def\accepted#1{\def\accepteddate{#1}}

\let\\\par\let\thelognumber\relax\let\thevolumenumber\relax
\let\thepapernumber\relax\let\thevolumeyear\relax\let\startpage\relax
\let\finishpage\relax\let\publishdate\relax\let\receiveddate\relax
\let\reviseddate\relax\let\accepteddate\relax\let\theasciititle\relax
\let\theasciiauthors\relax
\let\theasciiabstract\relax

\let\theasciiemail\relax

\ifplaintex
\font\logobig=cmssbx10 scaled 3836
\font\logomed=cmssbx10 scaled 2557
\else
\font\logobig=cmssbx10 scaled 4200
\font\logomed=cmssbx10 scaled 2800
\fi

\long\def\makeagttitle{   
\count0=\startpage
\agt\hfill      
\hbox to 45truept{\vbox to 0pt{\vglue -13truept{\logomed A\kern -.37em{\logobig 
T}\kern -.38em G}\vss}\hss}
\break
{\small Volume \thevolumenumber\ (\thevolumeyear)
\startpage--\finishpage\nl
Published: \publishdate}

\vglue .25truein

{\parskip=0pt\leftskip 0pt plus
1fil\def\\{\par\smallskip}{\Large\bf\thetitle}\par\medskip} \vglue
0.05truein

{\parskip=0pt\leftskip 0pt plus 1fil\def\\{\par}{\sc\theauthors}
\par\medskip}%
 
\vglue 0.03truein 

{\small\leftskip 25truept\rightskip 25truept{\bf Abstract}\stdspace\theabstract

{\bf AMS Classification}\stdspace\theprimaryclass
\ifx\thesecondaryclass\relax\else; \thesecondaryclass\fi\par
{\bf Keywords}\stdspace \thekeywords\par}\vglue 7truept

}   

\ifplaintex
\font\phead=cmsl9 scaled 950
\font\pnum=cmbx10 scaled 913
\font\pfoot=cmsl9 scaled 950
\headline{\vbox to 0pt{\vskip -4.5mm\line{\small\phead\ifnum
\count0=\startpage ISSN 1472-2739 (on-line) 1472-2747 (printed)
\hfill {\pnum\folio}\else\ifodd\count0\def\\{ }%
\ifx\theshorttitle\relax\thetitle\else\theshorttitle\fi\hfill{\pnum\folio}
\else\def\\{ and }{\pnum\folio}\hfill\ifx\theshortauthors\relax\theauthors
\else\theshortauthors\fi\fi\fi}\vss}}
\footline{\vbox to 0pt{\vglue 0mm\line{\small\pfoot\ifnum\count0=\startpage
\copyright\ \gtp\hfill\else
\agt, Volume \thevolumenumber\ (\thevolumeyear)\hfill\fi}\vss}}
\else
\font=cmsl9 scaled 1050
\font\lnum=cmbx10 
\font=cmsl9 scaled 1050
\makeatletter
\def\@oddhead{{\small\lhead\ifnum\count0=\startpage ISSN 1472-2739 
(on-line) 1472-2747 (printed)\hfill {\lnum\number\count0}\else\ifodd\count0
\def\\{ }\ifx\theshorttitle\relax \thetitle \else\theshorttitle\fi\hfill
{\lnum\number\count0}\else\def\\{ and }{\lnum\number\count0}
\hfill\ifx\theshortauthors\relax 
\theauthors\else\theshortauthors\fi\fi\fi}}\def\@evenhead{\@oddhead}
\def\@oddfoot{\small\lfoot\ifnum\count0=\startpage\copyright\ \gtp\hfill\else
\agt, Volume \thevolumenumber\ (\thevolumeyear)\hfill\fi}
\def\@evenfoot{\@oddfoot}
\makeatother
\fi
\let\maketitlepage\makeagttitle
\let\makeshorttitle\maketitlepage
\let\maketitle\maketitlepage

\newwrite\gtoutfile
\long\gdef\makeheadfile{  
{\def\\{, }\def\s{ }
\immediate\openout\gtoutfile head.xxx
\immediate\write\gtoutfile{Proxy-for: \ifx\theasciiauthors\relax
\theauthors\else\theasciiauthors\fi\s<\ifx\theasciiemail\relax\theemail\else\theasciiemail\fi>}
\immediate\write\gtoutfile{\noexpand\\}
\immediate\write\gtoutfile{Authors: \ifx\theasciiauthors\relax
\theauthors\else\theasciiauthors\fi}
{\def\\{ }\immediate\write\gtoutfile{Title: \ifx\theasciititle\relax
\thetitle\else\theasciititle\fi}}
\immediate\write\gtoutfile{Subj-class: GT or SG, GR etc}
\immediate\write\gtoutfile{MSC-class: \theprimaryclass\ifx\thesecondaryclass\relax\else, \thesecondaryclass\fi}
\immediate\write\gtoutfile{Journal-ref: Algebr. Geom. Topol. \thevolumenumber\s
(\thevolumeyear) \startpage-\finishpage}
\immediate\write\gtoutfile{Comments: Published by Algebraic and
Geometric Topology at}
\immediate\write\gtoutfile{\s\s\s  http://www.maths.warwick.ac.uk/agt/AGTVol\thevolumenumber/agt-\thevolumenumber-\thepapernumber.abs.html}
\immediate\write\gtoutfile{\noexpand\\}
\immediate\write\gtoutfile{}
\ifx\theasciiabstract\relax
\immediate\write\gtoutfile{\theabstract}\else
\immediate\write\gtoutfile{\theasciiabstract}\fi
\immediate\write\gtoutfile{}
\immediate\write\gtoutfile{\noexpand\\}
\immediate\write\gtoutfile{}
\immediate\closeout\gtoutfile}}  

\def\maketitlepage{\makeagttitle\makeheadfile}
\let\makeshorttitle\maketitlepage
\let\maketitle\maketitlepage